\documentclass[10pt]{amsart}

\newcommand{\proclaim}[2]{\medbreak {\bf #1}{\sl #2} \medbreak}

\def\cal{\mathcal}

\let\newpf\proof \let\proof\relax 
\newenvironment{pf}{\newpf[\proofname]}{\qed\endtrivlist}

\def\0{{\mathbf{0}}}

\def\r{\rho}
\def\tr{\tilde \rho}
\def\g{\gamma}
\def\tg{{\tilde\gamma}}

\def\be{\begin{equation}}
\def\ee{\end{equation}}

\def\d{{\underline d}}

\newcommand{\ntop}[2]{\genfrac{}{}{0pt}{1}{#1}{#2}}

\newtheorem{thm}{Theorem}[section]
\newtheorem{cor}[thm]{Corollary}

\newtheorem{lemma}[thm]{Lemma}

\newtheorem{prop}[thm]{Proposition}

\theoremstyle{remark}
\newtheorem{rem}{Remark}[section]

\numberwithin{equation}{section}

\def \bn {\hfill \\ \smallskip\noindent}

\theoremstyle{definition}

\def\proof{\bn {\bf Proof.} }

\newcommand{\dist}{\operatorname{dist}}

\newcommand{\inter}{\operatorname{int}}
\renewcommand{\mod}{\operatorname{mod}}

\newcommand{\id}{\operatorname{id}}

\newcommand{\conj}{\operatorname{conj}}

\newcommand{\Dil}{\operatorname{Dil}}

\newcommand{\CC}{{\cal C}}
\newcommand{\DD}{{\cal D}}
\newcommand{\EE}{{\cal E}}

\newcommand{\II}{{\cal I}}
\newcommand{\FF}{{\cal F}}

\newcommand{\KK}{{\cal K}}

\newcommand{\QQ}{{\cal Q}}
\newcommand{\RR}{{\cal R}}
\newcommand{\TT}{{\cal T}}

\newcommand{\UU}{{\cal U}}
\newcommand{\VV}{{\cal V}}
\newcommand{\WW}{{\cal W}}
\newcommand{\XX}{{\cal X}}
\newcommand{\YY}{{\cal Y}}

\newcommand{\hh}{{\mathbf h}}

\newcommand{\C}{{\mathbb C}}
\newcommand{\D}{{\mathbb D}}

\newcommand{\N}{{\mathbb N}}

\newcommand{\R}{{\mathbb R}}

\newcommand{\Z}{{\mathbb Z}}

\def\B0{{\bold{0}}}

\catcode`\@=12

\def\Empty{}
\newcommand\oplabel[1]{
  \def\OpArg{#1} \ifx \OpArg\Empty {} \else
        \label{#1}
  \fi}
                
\newcommand{\comm}[1]{}
\newcommand{\comment}[1]{}

\begin{document}

\title[Statistical properties in the quadratic family]
{Statistical properties of unimodal maps: the quadratic family}

\author{Artur Avila and Carlos Gustavo Moreira}

\address{
Coll\`ege de France -- 3 Rue d'Ulm \\
75005 Paris -- France.
}
\email{avila@impa.br}

\address{
IMPA -- Estr. D. Castorina 110 \\
22460-320 Rio de Janeiro -- Brazil.
}
\email{gugu@impa.br}

\thanks{Partially supported by Faperj and CNPq, Brazil.}

\date{\today}

\begin{abstract}

We prove that almost every non-regular
real quadratic map is Collet-Eckmann and has
polynomial recurrence of the critical orbit (proving a
conjecture by Sinai).
It follows that typical quadratic maps have excellent ergodic
properties, as exponential decay of correlations
(Keller and Nowicki, Young) and stochastic stability in the
strong sense (Baladi and Viana).
This is an important step to get the same results for
more general families of unimodal maps.

\end{abstract}

\setcounter{tocdepth}{1}

\maketitle

\tableofcontents

\section*{Introduction}

Here we consider the quadratic family, $f_a=a-x^2$, where $-1/4 \leq a \leq
2$ is the parameter, and we analyze its dynamics in the invariant interval.

The quadratic family has been one of the most studied dynamical systems in
the last decades.  It is one of the most basic examples and exhibits a very
rich behavior.  It was also studied through many different techniques.
Here we are interested in describing the dynamics of a typical quadratic map
from the statistical point of view.

\subsection{The probabilistic point of view in dynamics}

In the last decade Palis \cite {Pa} described a general program for
(dissipative) dynamical systems in any dimension.  In short, he asks to show
that `typical' dynamical systems can be modeled stochastically in a
robust way.  More precisely, one should show that such typical system
can be described by finitely many attractors, each of them supporting
a (ergodic) physical measure: time averages of Lebesgue almost every orbit
should converge to spatial averages according to one of the
physical measures.  The description should be robust
under (sufficiently) random perturbations of the system: one asks for
stochastic stability.

Moreover,
a typical dynamical system was to be understood in the Kolmogorov
sense: as a set of full measure in generic parametrized families.

Besides the questions posed by this conjecture, much more can be asked about
the statistical description of the long term behavior of a typical system.
For instance, the definition of physical measure is related to the
validity of the Law of Large Numbers.  Are other theorems still valid,
like the Central Limit or Large Deviation theorems?  Those questions
are usually related to the rates of mixing of the physical measure.

\subsection{The richness of the quadratic family}

While we seem still very far away of any description of dynamics of typical
dynamical systems (even in one-dimension),
the quadratic family has been a remarkable exception.  Let us describe
briefly some results which show the richness of the quadratic family
from the probabilistic point of view.

The initial step in this direction was the work of Jakobson \cite {J},
where it was shown that for a positive measure set of parameters
the behavior is stochastic, more precisely, there is an
absolutely continuous invariant measure
(the physical measure) with positive Lyapunov exponent: for Lebesgue
almost every $x$, $|Df^n(x)|$ grows exponentially fast.
On the other hand, it was later shown by
Lyubich \cite {puzzle} and Graczyk-Swiatek \cite {GS} that regular
parameters (with a periodic hyperbolic attractor) are (open and) dense.
While stochastic parameters are predominantly expanding
(in particular have sensitive dependence to initial conditions), regular
parameters are deterministic (given by the periodic attractor).
So at least two kinds of very distinct observable behavior are present
in the quadratic family, and they alternate in a complicate way.

It was later shown that stochastic behavior could be concluded from enough
expansion along the orbit of the critical value: the Collet-Eckmann
condition, exponential growth of $|Df^n(f(0))|$, was enough to conclude
a positive Lyapunov exponent of the system.
A different approach to Jakobson's Theorem in \cite {BC1} and \cite {BC2}
focused specifically on this property: the set of
Collet-Eckmann maps has positive
measure.  After these initial works, many others studied such parameters
(sometimes with extra assumptions), obtaining refined information
of the dynamics of CE maps, particularly information about
exponential decay of correlations\footnote
{CE quadratic maps are not always mixing and finite periodicity can
appear in a robust way.  This phenomena is related to the map being
renormalizable, and this is the only obstruction:
the system is exponentially mixing after renormalization.}
(Keller and Nowicki
in \cite {KN} and Young in \cite {Y2}), and stochastic stability (Baladi
and Viana in \cite {BV}).  The dynamical systems considered in those papers
have generally been shown to have excellent statistical
description\footnote
{It is now known that weaker expansion than Collet-Eckmann
is enough to obtain stochastic behavior for quadratic maps,
on the other hand, exponential decay of correlations is actually
equivalent to the CE condition \cite {NS}, and all current results on
stochastic stability use the Collet-Eckmann condition.}.

Many of those results also generalized to more general
families and sometimes to higher dimensions, as in the case of
H\'enon maps \cite {BC2}.

The main motivation behind this strong effort to understand the class of
CE maps was certainly the fact that such a class was known to have
positive measure.  It was known however that very different
(sometimes wild) behavior coexisted.  For instance, it was shown the
existence of quadratic maps without a physical measure or quadratic maps
with a physical measure concentrated on a {\it repelling} hyperbolic
fixed point (\cite {Jo}, \cite {HK}).  It remained to see if wild behavior
was observable.

In a big project in the last decade, Lyubich \cite {parapuzzle} together with
Martens and Nowicki \cite {MN} showed that almost all parameters have
physical measures: more precisely, besides regular and stochastic behavior,
only one more behavior could (possibly) happen with positive measure, namely
infinitely renormalizable maps (which always have a uniquely ergodic
physical measure).  Later Lyubich in \cite {regular} showed that infinitely
renormalizable parameters have measure zero, thus establishing the
celebrated {\it regular or stochastic} dichotomy.  This further advancement
in the comprehension of the nature of the statistical behavior of typical
quadratic maps is remarkably linked to the progress obtained by
Lyubich on the answer of the Feigenbaum conjectures \cite {universe}.

\subsection{Statements of the results}

In this work we describe the asymptotic behavior of the critical orbit. 
Our first result is an estimate of hyperbolicity:
\proclaim{Theorem A.}
{
Almost every non-regular real quadratic map satisfies the Collet-Eckmann
condition:
$$
\liminf_{n \to \infty} \frac {\ln(|Df^n(f(0))|)} {n}>0.
$$
}

The second is an estimate on the recurrence of the critical point.  For
regular maps, the critical point is non-recurrent (it actually converges to
the periodic attractor).  Among non-regular maps, however, the recurrence  
occurs at a precise rate which we estimate:
\proclaim{Theorem B.}
{
Almost every non-regular real quadratic map
has polynomial recurrence of the critical orbit with exponent $1$:
$$
\limsup_{n \to \infty} \frac {-\ln(|f^n(0)|)} {\ln(n)}=1.
$$
In other words, the set of $n$ such that $|f^n(0)|<n^{-\g}$ is finite
if $\g>1$ and infinite if $\g<1$.
}

As far as we know, this is the first proof of polynomial estimates for
the recurrence of the critical orbit valid for a positive measure set of
non-hyperbolic parameters (although subexponential
estimates were known before).  This
also answers a long standing conjecture of Sinai.

Theorems A and B show that a typical non regular quadratic maps have enough
good properties to conclude the results on
exponential decay of correlations (which can be used to prove Central Limit
and Large Deviation theorems) and stochastic stability in the sense of $L^1$
convergence of the densities (of stationary measures of perturbed systems).
Many other properties also follow, like existence
of a spectral gap in \cite {KN} and the recent results on
almost sure (stretched exponential) rates of convergence
to equilibrium in \cite {BaBeM}.  In particular, this answers positively
Palis conjecture for the quadratic family.

\subsection{Unimodal maps}

Another reason to deal with the quadratic family is that it seems to open
the doors to the understanding of unimodal maps.  Its universal behavior was
first realized in the topological sense, with Milnor-Thurston theory.  The
Feigenbaum-Coullet-Tresser observations indicated a geometric universality
\cite {universe}.

A first result in the understanding of measure-theoretical universality was
the work of Avila, Lyubich and de Melo
\cite {ALM}, where it was shown how to relate metrically
the parameter spaces of non-trivial
analytic families of unimodal maps to the parameter space of the quadratic
family.  This was proposed as a method
to relate observable dynamics in the
quadratic family to observable dynamics of general analytic families
of unimodal maps.  In that work the method is used successfully to extend the
regular or stochastic dichotomy to this broader context.

We are also able to adapt those methods to our setting.  The techniques
developed here and the methods of \cite {ALM} are the main tools
used in \cite {AM} to obtain the main results of this paper
(except the exact value of the polynomial recurrence)
for non-trivial real analytic families of unimodal maps
(with negative Schwarzian derivative and quadratic
critical point).  This is a rather general set of families,
as trivial families form a set of infinite codimension.  For a different
approach (still based on \cite {ALM})
which does not use negative Schwarzian derivative and obtains
the exponent $1$ for the polynomial recurrence, see \cite {A}.

In \cite {AM} we also prove a version of Palis conjecture
in the smooth setting.  There is a residual set of $k$-parameter
$C^3$ (for the equivalent $C^2$ result, see \cite {A})
families of unimodal maps with negative Schwarzian derivative
such that almost every parameter is either regular or Collet-Eckmann
with subexponential bounds for the recurrence of the critical point.

{\bf Acknowledgements:}
We thank Viviane Baladi, Mikhail Lyubich, Marcelo Viana, and Jean-Christophe
Yoccoz for helpful discussions.  We are grateful to
Juan Rivera-Letelier for listening to a first version, and for valuable
discussions on the Phase-Parameter relation, which led to the use of the
gape interval in this work.  We would like to thank the anonymous referee
for his suggestions concerning the presentation of this paper.

\section{General definitions}

\subsection{Maps of the interval}

Let $f:I \to I$ be a $C^1$ map defined on some interval $I \subset \R$.
The {\it orbit} of a point $p \in I$ is the sequence
$\{f^k(p)\}_{k=0}^\infty$.
We say that $p$ is {\it recurrent} if there exists a
subsequence $n_k \to \infty$
such that $\lim f^{n_k}(p)=p$.

We say that $p$ is a {\it periodic point of period $n$} of $f$ if
$f^n(p)=p$, and $n \geq 1$ is minimal with this property.
In this case we say that $p$ is
{\it hyperbolic} if $|Df^n(p)|$ is not $0$ or $1$.
Hyperbolic periodic orbits are
{\it attracting} or {\it repelling} according to
$|Df^n(p)|<1$ or $|Df^n(p)|>1$.

We will often consider the restriction of iterates $f^n$ to
intervals $T \subset I$, such that $f^n|_T$ is a diffeomorphism.  In
this case we will be interested on the {\it distortion} of $f^n|_T$,
$$
\dist(f^n|_T)=\frac {\sup_T |Df^n|} {\inf_T |Df^n|}.
$$
This is always a number bigger than or equal to $1$, we will say that it is
small if it is close to $1$.

\subsection{Trees}

We let $\Omega$ denote the set of finite sequences of non-zero integers
(including the empty sequence).  Let $\Omega_0$ denote $\Omega$ without the
empty sequence.  For $\d \in \Omega$, $\d=(j_1,...,j_m)$, we let $|\d|=m$
denote its length.

We denote $\sigma^+:\Omega_0 \to \Omega$ by
$\sigma^+(j_1,...,j_m)=(j_1,...,j_{m-1})$ and
$\sigma^-:\Omega_0 \to \Omega$ by
$\sigma^-(j_1,...,j_m)=(j_2,...,j_m)$.

For the purposes of this paper, one should view $\Omega$ as a (directed)
tree with root $\d=\emptyset$ and edges connecting $\sigma^+(\d)$ to $\d$
for each $\d \in \Omega_0$.  We will use $\Omega$ to label objects which are
organized in a similar tree structure (for instance, certain families of
intervals ordered by inclusion).

\subsection{Growth of functions}

Let $f:\N \to \R^+$ be a function.  We say that $f$ grows {\it at least
exponentially} if there exists $\alpha>0$ such that $f(n)>e^{\alpha n}$ for
all $n$ sufficiently big.  We
say that $f$ grows {\it at least polynomially} if there exists $\alpha>0$
such that $f(n)>n^\alpha$ for all $n$ sufficiently big.

The standard {\it torrential} function $T$ is defined recursively by
$T(1)=1$, $T(n+1)=2^{T(n)}$.  We say that $f$ grows {\it at
least torrentially} if there exists $k>0$ such that $f(n)>T(n-k)$ for every
$n$ sufficiently big.  We will say
that $f$ grows {\it torrentially} if there exists
$k>0$ such that $T(n-k)<f(n)<T(n+k)$ for every $n$ sufficiently big.

Torrential growth can be detected from recurrent estimates easily.
A sufficient condition for a function $f$
to grow at least torrentially is an estimate as
$$
f(n+1)>e^{f(n)^\alpha}
$$
for some $\alpha>0$.  Torrential growth is implied by an estimate as
$$
e^{f(n)^\alpha}<f(n+1)<e^{f(n)^\beta}
$$
with $0<\alpha<\beta$.

We will also say that $f$ decreases at least exponentially (respectively
torrentially) if $1/f$ grows at least exponentially (respectively
torrentially).

\comm{
Let $X$ be a class of functions $g:\N \to \R$ such that $\lim_{n \to \infty}
g(n)=\infty$.  We say that a function $f:\N \to \R$ grows with rate at least
$X$ if there exists a function $g \in X$ such that $f(n) \geq g(n)$ for $n$
sufficiently big.  We say
that it grows at rate $X$ if there are
$g_1,g_2 \in X$ such that $g_1(n)
\leq f(n) \leq g_2(n)$ for $n$ sufficiently big.  We say that $f$ decreases
with rate (at least) $X$ if
$1/f$ grows at rate (at least) $X$.

Standard classes are the following.  Linear for linear functions with
positive slope.  Polynomial for functions $g(n)=n^k,k>0$.  Exponential for
functions $g(n)=e^{k n},k>0$.

The standard {\it torrential} function $T$ is defined recursively by
$T(1)=1$, $T(n+1)=2^{T(n)}$.  The torrential class is the set of functions
$g(n)=T(\max\{n+k,1\}), k \in \Z$.

Torrential growth can be detected from recurrent estimates easily.
A sufficient condition for a function which is unbounded from above
to grow at least torrentially is an estimate as
$$
f(n+1)>e^{f(n)^a}
$$
for some $a>0$.  Torrential growth is implied by an estimate as
$$
e^{f(n)^a}<f(n+1)<e^{f(n)^b}
$$
with $0<a<b$.
}

\subsection{Quasisymmetric maps}

Let $k \geq 1$ be given.
We say that a homeomorphism $f:\R \to \R$ is {\it quasisymmetric}
with constant $k$ if for all $h > 0$
$$
\frac {1} {k} \leq \frac {f(x+h)-f(x)} {f(x)-f(x-h)} \leq k.
$$

The space of quasisymmetric maps is a group under composition, and the set
of quasisymmetric maps with constant $k$ preserving a given interval is
compact in the uniform topology of compact subsets of $\R$.  It also follows
that quasisymmetric maps are H\"older.

To describe further the properties of quasisymmetric maps, we need the
concept of quasiconformal maps and dilatation
so we just mention a result of Ahlfors-Beurling
which connects both concepts: any quasisymmetric map extends
to a quasiconformal real-symmetric map of
$\C$ and, conversely, the
restriction of a quasiconformal real-symmetric map of $\C$ to $\R$ is
quasisymmetric.  Furthermore, it is possible to work out upper bounds on the
dilatation (of an optimal extension) depending only on $k$ and conversely.

The constant $k$ is awkward to work with: the inverse of a quasisymmetric
map with constant $k$ may have a larger constant.  We will therefore work
with a less standard constant: we will say that $h$ is
$\g$-quasisymmetric ($\g$-qs) if $h$ admits a quasiconformal symmetric
extension to $\C$ with dilatation bounded by $\g$.  This definition
behaves much better: if $h_1$ is $\g_1$-qs and $h_2$ is
$\g_2$-qs then $h_2 \circ h_1$ is $\g_2 \g_1$-qs.

If $X \subset \R$ and $h:X \to \R$ has a $\g$-quasisymmetric extension to
$\R$ we will also say that $h$ is $\g$-qs.

Let $QS(\g)$ be the set of $\g$-qs maps of $\R$.

\section{Real quadratic maps}

If $a \in \C$ we let $f_a:\C \to \C$ denote the (complex)
quadratic map $a-z^2$.  For real parameters in the range
$-1/4 \leq a \leq 2$, there exists an interval
$I_a=[\beta,-\beta]$ with
$$
\beta=\frac {-1-\sqrt {1+4 a}} {2}
$$
such that $f_a(I_a) \subset I_a$ and $f_a(\partial
I_a) \subset \partial I_a$.  For such values of the parameter $a$,
the map $f=f_a|_{I_a}$ is unimodal, that is, it is a self map of
$I_a$ with a unique turning point.  To simplify the notation, we will
usually drop the dependence on the parameter and let $I=I_a$.



\subsection{The combinatorics of unimodal maps}

In this subsection we fix a real quadratic map $f$ and define some objects
related to it.

\subsubsection{Return maps}

Given an interval $T \subset I$ we define the {\it first return map}
$R_T:X \to T$ where $X \subset T$ is the set of points
$x$ such that there exists $n>0$ with
$f^n(x) \in T$, and $R_T(x)=f^n(x)$ for the minimal $n$ with this property.

\subsubsection{Nice intervals}

An interval $T$ is {\it nice}
if it is symmetric around $0$ and the iterates of
$\partial T$ never intersect $\inter T$.  Given a nice interval $T$
we notice that the domain of the
first return map $R_T$ decomposes in a union of intervals
$T^j$, indexed by integer numbers (if there are only finitely many
intervals, some indexes will be corresponded to the empty set).
If $0$ belongs to the domain of $R_T$, we say that $T$ is {\it proper}.
In this case we reserve the index $0$ to denote the
component of the critical point: $0 \in T^0$.

If $T$ is nice, it follows that for all $j \in \Z$, $R_T(\partial T^j)
\subset \partial T$.  In particular, $R_T|_{T^j}$ is
a diffeomorphism onto $T$ unless $0 \in T^j$ (and in particular $j=0$ and
$T$ is proper).  If $T$ is proper, $R_T|_{T^0}$ is symmetric (even) with a
unique critical point $0$.  As a
consequence, $T^0$ is also a nice interval.

If $R_T(0) \in T^0$, we say that $R_T$ is {\it central}.

If $T$ is a proper interval then both $R_T$ and $R_{T^0}$ are defined, and
we say that $R_{T^0}$ is the generalized renormalization of $R_T$.

\subsubsection{Landing maps}

Given a proper interval $T$ we define the {\it landing map} $L_T:X \to T^0$
where $X \subset T$ is the set of points $x$ such that there exists
$n \geq 0$ with
$f^n(x) \in T^0$, and $L_T(x)=f^n(x)$ for the minimal
$n$ with this property.  We notice that $L_T|_{T^0}=\id$.

\subsubsection{Trees}

We will use $\Omega$ to label iterations of non-central branches of $R_T$,
as well as their domains.
If $\d \in \Omega$, we define
$T^\d$ inductively in the following way.  We let $T^\d=T$
if $\d$ is empty and if $\d=(j_1,...,j_m)$ we let
$T^\d=(R_T|_{T^{j_1}})^{-1}(T^{\sigma^-(\d)})$.

We denote $R^\d_T=R_T^{|\d|}|_{T^\d}$
which is always a diffeomorphism onto $T$.

Notice that the family of intervals $T^\d$ is organized by inclusion in
the same way as $\Omega$ is organized by (right side) truncation
(the previously introduced tree structure).

If $T$ is a proper interval, the first return map to $T$ naturally relates
to the first landing to $T^0$.  Indeed, denoting $C^\d=(R^\d_T)^{-1}(T^0)$,
the domain of the first landing map $L_T$ is easily seen to coincide with
the union of the $C^\d$, and furthermore $L_T|_{C^\d}=R^\d_T$.

Notice that this allows us to relate $R_T$ and $R_{T^0}$ since
$R_{T^0}=L_T \circ R_T$.

\subsubsection{Renormalization}

We say that $f$ is {\it renormalizable} if there is an interval $0 \in T$
and $m>1$ such that $f^m(T) \subset T$ and
$f^j(\inter T) \cap \inter T=\emptyset$ for $1 \leq j<m$.
The maximal such interval is called the {\it renormalization interval of
period $m$}, it has the property that $f^m(\partial T) \subset \partial T$.

The set of renormalization periods of $f$ gives an increasing (possibly
empty) sequence of numbers $m_i$, $i=1,2,...$,
each related to a unique renormalization
interval $T^{(i)}$ which form a nested sequence of intervals.  We include
$m_0=1$, $T^{(0)}=I$ in the sequence to simplify the notation.

We say that $f$ is {\it finitely renormalizable} if there is a smallest
renormalization interval $T^{(k)}$.  We say that
$f \in \FF$ if $f$ is finitely
renormalizable and $0$ is recurrent but not periodic.  We let $\FF_k$ denote
the set of maps $f$ in $\FF$ which are exactly $k$ times renormalizable.

\subsubsection{Principal nest}

Let $\Delta_k$ denote the set of all maps $f$
which have (at least) $k$ renormalizations and which have an orientation
reversing non-attracting periodic point of period $m_k$ which we denote
$p_k$ (that is, $p_k$ is the fixed point of
$f^{m_k}|_{T^{(k)}}$ with $Df^{m_k}(p_k) \leq -1$).
For $f \in \Delta_k$, we denote $T^{(k)}_0=[-p_k,p_k]$.
We define by induction a (possibly
finite) sequence $T^{(k)}_i$, such that $T^{(k)}_{i+1}$ is the component of
the domain of $R_{T^{(k)}_i}$ containing $0$.  If this sequence is
infinite, then either it converges to a point or to an interval.

If $\cap_i T^{(k)}_i$ is a point, then
$f$ has a recurrent critical point which is not
periodic, and it is possible to show that $f$ is not $k+1$ times
renormalizable.  Obviously in this case we have
$f \in \FF_k$, and all maps in $\FF_k$ are obtained in this way:
if $\cap_i T^{(k)}_i$ is an
interval, it is possible to show that $f$ is $k+1$ times renormalizable.

We can of course write $\FF$ as a disjoint union
$\cup_{i=0}^\infty \FF_i$.  For a map
$f \in \FF_k$ we refer to the sequence
$\{T^{(k)}_i\}_{i=1}^\infty$ as the {\it principal nest}.

It is important to notice that the domain of the first return map
to $T^{(k)}_i$ is always dense in $T^{(k)}_i$.  Moreover, the next result
shows that, outside a very special case,
the return map has a hyperbolic structure.

\begin{lemma} \label {hyperbol}

Assume $T^{(k)}_i$ does not have a non-hyperbolic periodic orbit in its
boundary.  For all $T^{(k)}_i$ there exists $C>0$, $\lambda>1$
such that if $x,f(x),...,f^{n-1}(x)$ do not belong to
$T^{(k)}_i$ then $|Df^n(x)|>C \lambda^n$.

\end{lemma}

This lemma is a simple consequence of a general theorem of Guckenheimer
on hyperbolicity of maps of the interval without critical points and
non-hyperbolic periodic orbits (Guckenheimer considers unimodal maps with
negative Schwarzian derivative, so this applies directly
to the case of quadratic maps, the general case is also true by
Ma\~n\'e's Theorem, see \cite {MvS}).  Notice that the existence of a
non-hyperbolic periodic orbit in the boundary of $T^{(k)}_i$ depends on a
very special combinatorial setting, in particular, all $T^{(k)}_j$ must
coincide (with $[-p_k,p_k]$), and the $k$-th renormalization of $f$ is in
fact renormalizable of period $2$.

By Lemma \ref {hyperbol}, the maximal
invariant of $f|_{I \setminus T^{(k)}_i}$ is an expanding set,
which admits a Markov partition (since $\partial T^{(k)}_i$ is
preperiodic, see also the proof of Lemma \ref {expan}): it is easy to
see that it is indeed a Cantor
set\footnote {Dynamically defined Cantor sets with such properties are
usually called {\it regular Cantor sets}.}
(except if $i=0$ or in the special
period $2$ renormalization case just described).
It follows that the geometry of this Cantor set is well behaved:
for instance, its image by any quasisymmetric map
has zero Lebesgue measure.

In particular, one sees that the domain of the
first return map to $T^{(k)}_i$ has infinitely many components (except in
the special case above or if $i=0$) and that its complement has well behaved
geometry.


\subsubsection{Lyubich's Regular or Stochastic dichotomy}

A map $f \in \FF_k$ is called {\it simple} if the principal nest has
only finitely many central returns, that is, there are
only finitely many $i$ such that $R|_{T^{(k)}_i}$ is central.
Such maps have many good features, in particular, they are stochastic
(this is a consequence of \cite {MN} and \cite {attractors}).

In \cite {parapuzzle}, it was proved that almost every
quadratic map is either regular or simple or infinitely renormalizable.  It
was then shown in \cite {regular} that infinitely renormalizable maps have
zero Lebesgue measure, which establishes the Regular or Stochastic
dichotomy.

Due to Lyubich's results, we can completely forget about infinitely
renormalizable maps, we just have to prove the claimed estimates for almost
every simple map.

During our discussion, for notational reasons, we will
fix a renormalization level $\kappa$, that is, we will
only analyze maps in $\Delta_\kappa$.
This allows us to fix some convenient notation: given $g \in \Delta_\kappa$
we define $I_i[g]=T^{(\kappa)}_i[g]$, so that $\{I_i[g]\}$
is a sequence of intervals (possibly finite).
We use the notation $R_i[g]=R_{I_i[g]}$, $L_i[g]=L_{I_i[g]}$ and so on (so
that the domain of $R_i[g]$ is $\cup I^j_i[g]$ and the domain of $L_i[g]$ is
$\cup C^\d_i[g]$).  When doing phase analysis (working with fixed $f$)
we usually drop the
dependence on the map and write $R_i$ for $R_i[f]$.

(Notice that, once we fix the renormalization level $\kappa$, for $g \in
\Delta_\kappa$, the notation $I_i[g]$ stands for $T^{(\kappa)}_i[g]$,
even if $g$ is
more than $\kappa$ times renormalizable.)

\subsubsection{Strategy}

To motivate our next steps, let us describe the general strategy behind the
proofs of Theorems A and B.

(1)\, We consider a certain set of non-regular parameters of full measure
and describe (in a probabilistic way) the dynamics of
the principal nest.  This is our phase analysis.

(2)\, From time to time, we transfer the information from the
phase space to the parameter, following the description of the parapuzzle
nest which we will make in the next subsection.  The rules for this
correspondence are referred as Phase-Parameter relation (which is based on
the work of Lyubich on complex dynamics of the quadratic family).

(3)\, This correspondence will
allow us to exclude parameters whose critical orbit behaves badly (from
the probabilistic point of view) at infinitely many levels of the principal
nest.  The phase analysis coupled with the
Phase-Parameter relation will assure us that the remaining parameters have
still full measure.

(4)\, We restart the phase analysis
for the remaining parameters with extra information.

After many iterations of this procedure we will have enough information to
tackle the problems of hyperbolicity and recurrence.

We will first tackle the problem of describing the Phase-Parameter relation,
and we will delay all statistical arguments until \S \ref {measure}.

A larger outline of this strategy, including the motivation and
organization of the statistical analysis, will appear in \cite {AM5}.

\subsection{Parameter partition}

Part of our work is to transfer information from the phase space of some
map $f \in \FF$ to a neighborhood of $f$ in the parameter space.
This is done in the following way.  We consider the
first landing map $L_i$: the complement of the domain of $L_i$ is a
hyperbolic Cantor set $K_i=I_i \setminus \cup C^\d_i$.
This Cantor set persists in a small parameter
neighborhood $J_i$ of $f$, changing in a continuous way.
Thus, loosely speaking, the domain of $L_i$ induces a persistent partition
of the interval $I_i$.

Along $J_i$, the first landing map is topologically the same (in a way
that will be clear soon).  However the critical value
$R_i[g](0)$ moves relative
to the partition (when $g$ moves in $J_i$).  This allows us to partition
the parameter piece $J_i$ in smaller pieces, each
corresponding to a region where $R_i(0)$ belongs to some fixed component of
the domain of the first landing map.

\begin{thm} [Topological Phase-Parameter relation]

Let $f \in \FF_\kappa$.
There is a sequence $\{J_i\}_{i \in \N}$ of nested
parameter intervals (the {\it principal parapuzzle nest} of $f$)
with the following properties.

\begin{enumerate}

\item $J_i$ is the maximal interval containing $f$ such that for all $g \in
J_i$ the interval
$I_{i+1}[g]=T^{(\kappa)}_{i+1}[g]$ is defined and changes in a
continuous way.  (Since the first return map to $R_i[g]$ has a central
domain, the landing map
$L_i[g]:\cup C^\d_i[g] \to I_i[g]$ is defined.)

\item $L_i[g]$ is topologically the same along $J_i$: there exists
homeomorphisms $H_i[g]:I_i \to I_i[g]$, such that
$H_i[g](C^\d_i)=C^\d_i[g]$.
The maps $H_i[g]$ may be chosen to change continuously.

\item There exists a homeomorphism $\Xi_i:I_i \to J_i$ such that
$\Xi_i(C^\d_i)$ is
the set of $g$ such that $R_i[g](0)$ belongs to $C^\d_i[g]$.

\end{enumerate}

\end{thm}

The homeomorphisms $H_i$ and $\Xi_i$ are not uniquely defined, it is easy to
see that we can modify them inside each $C^\d_i$ window keeping the above
properties.  However, $H_i$ and $\Xi_i$
are well defined maps if restricted to $K_i$.

This fairly standard phase-parameter result can be proved in many different
ways.  The most elementary proof is
probably to use the monotonicity of the
quadratic family to deduce the Topological Phase-Parameter relation from
Milnor-Thurston's kneading theory by purely combinatorial arguments.
Another approach is to use Douady-Hubbard's description of the combinatorics
of the Mandelbrot set (restricted to the real line) as does Lyubich
in \cite {parapuzzle} (see also \cite {AM4} for a more general case).

With this result we can define for any $f \in \FF_\kappa$ intervals
$J^j_i=\Xi_i(I^j_i)$ and $J^\d_i=\Xi_i(I^\d_i)$.  From the description we
gave it immediately follows that two intervals $J_{i_1}[f]$ and $J_{i_2}[g]$
associated to maps $f$ and $g$ are either disjoint or nested,
and the same happens for
intervals $J^j_i$ or $J^\d_i$.  Notice that if $g \in \Xi_i(C^\d_i) \cap
\FF_\kappa$ then $\Xi_i(C^\d_i)=J_{i+1}[g]$.

We will concentrate on the analysis of the regularity of $\Xi_i$ for the
special class of simple maps $f$: one of the good properties of the class of
simple maps is better control of the phase-parameter relation.
Even for simple maps, however, the regularity of $\Xi_i$ is not great:
there is too much dynamical information contained in it.  A solution to this
problem is to forget some dynamical information.

\subsubsection{Gape interval} \label {gap}

If $i>1$, we define the {\it gape interval} $\tilde I_{i+1}$ as follows.

We have that
$R_i|_{I_{i+1}}=L_{i-1} \circ R_{i-1}=R^{\d}_{i-1} \circ R_{i-1}$
for some $\d$, so that
$I_{i+1}=(R_{i-1}|_{I_i})^{-1}(C^{\d}_{i-1})$.  We define the
gape interval
$\tilde I_{i+1}=(R_{i-1}|_{I_i})^{-1}(I^{\d}_{i-1})$.

Notice that $I_{i+1} \subset \tilde I_{i+1} \subset I_i$.  Furthermore,
for each $I^j_i$, the gape interval
$\tilde I_{i+1}$ either contains or is disjoint from $I^j_i$.

\subsubsection{The Phase-Parameter relation}

As we discussed before, the dynamical information contained in $\Xi_i$
is entirely given by $\Xi_i|_{K_i}$: a map obtained by $\Xi_i$ by
modification inside a $C^\d_i$ window has still the same properties. 
Therefore it makes sense to ask about the regularity of $\Xi_i|_{K_i}$.  As
we anticipated before we must erase some information to obtain good results. 

Let $f \in \FF_\kappa$ and let $\tau_i$ be such that
$R_i(0) \in I^{\tau_i}_i$.  We define two Cantor sets,
$K^\tau_i=K_i \cap I^{\tau_i}_i$ which
contains refined information restricted to the $I^{\tau_i}_i$ window and
$\tilde K_i=I_i \setminus (\cup I^j_i \cup \tilde I_{i+1})$, which
contains global information, at the cost of erasing information inside each
$I^j_i$ window and in $\tilde I_{i+1}$.

\begin{thm} [Phase-Parameter relation]

Let $f$ be a simple map.  For all $\g>1$ there exists $i_0$ such that
for all $i>i_0$ we have

\begin{description}

\item[PhPa1] $\Xi_i|_{K^\tau_i}$ is $\g$-qs,

\item[PhPa2] $\Xi_i|_{\tilde K_i}$ is $\g$-qs,

\item[PhPh1] $H_i[g]|_{K_i}$ is $\g$-qs if $g \in J^{\tau_i}_i$,

\item[PhPh2] the map $H_i[g]|_{\tilde K_i}$ is
$\g$-qs if $g \in J_i$.

\end{description}

\end{thm}

The Phase-Parameter relation follows from the work of
Lyubich \cite {parapuzzle}, where a general method based on the
theory of holomorphic motions was introduced to deal with this kind of
problem.  A sketch of the derivation
of the specific statement of the Phase-Parameter relation from the general
method of Lyubich is given in the
Appendix A.  The reader can find full details (in a more general context
than quadratic maps) in \cite {AM4}.

\comm{
The proof of this Theorem will be given in \S \ref {pp}, but involves mainly
complex dynamics, so we now turn to discuss the complex
analogous of the principal nest description of unimodal maps.
}

\begin{rem}

One of the main reasons why the present work is
restricted to the quadratic family is related to
the Topological Phase-Parameter relation and the
Phase-Parameter relation.  The work of Lyubich uses specifics of the
quadratic family, specially the fact that it is a full family of
quadratic-like maps, and several arguments involved have indeed
a global nature (using for instance the combinatorial theory of
the Mandelbrot set).  Thus we are only able to conclude
the Phase-Parameter relation in this restricted setting.

However, the statistical analysis involved in
the proofs of Theorem A and B in this work is valid in much more generality.
Our arguments suffice (without any changes)
for any one-parameter analytic family of unimodal maps $f_\lambda$
with the following properties:

\noindent(1)\,
For every $\lambda$, $f_\lambda$ has a quadratic critical point and negative
Schwarzian derivative\footnote {More generally it is enough to ask that the
first return map to a sufficiently small nice interval has negative
Schwarzian derivative.},

\noindent(2)\,
For almost every non-regular parameter $\lambda$, $f_\lambda$ has all
periodic orbits repelling (so that Lemma \ref {hyperbol} holds), is
conjugate to a quadratic simple map, and the Topological
Phase-Parameter relation\footnote
{Actually one only needs the Topological Phase-Parameter relation
to be valid for all deep enough levels of the principal nest.}
and the Phase-Parameter relation\footnote
{In \cite {AM} it is shown how to work around this condition for
most families satisfying condition (1).  The results
obtained are weaker though, and the statistical analysis is slightly
harder.}
are valid at $\lambda$.



The assumption of a quadratic critical point is probably the hardest
to remove at this point, so our analysis does not apply, say,
for the families $a-x^{2n}$, $n>1$.  It is worth to point out that most of
the arguments developed in this paper go through for higher criticality. 
The key missing links are in the starting points of this paper:
zero Lebesgue measure of infinitely
renormalizable parameters and of finitely renormalizable
parameters without exponential decay of geometry (in the sense of \cite
{attractors}), and growth of moduli of parapuzzle annuli
(in the sense of \cite {parapuzzle}) for almost every parameter.

\end{rem}

\section{Measure and capacities} \label {measure}

\subsection{Quasisymmetric maps}

If $X \subset \R$ is measurable, let us denote $|X|$ its Lebesgue measure.
Let us explicit the metric properties of $\g$-qs maps we will use.

To each $\g$, there exists a constant $k \geq 1$ such that for all
$f \in QS(\g)$, for all $J \subset I$ intervals,
$$
\frac {1} {k} \left ( \frac {|J|} {|I|} \right )^k \leq \frac {|f(J)|}
{|f(I)|} \leq
\left ( \frac {k|J|} {|I|} \right )^{1/k}.
$$

Furthermore $\lim_{\g \to 1} k(\g)=1$.
So for each $\epsilon>0$ there exists $\g>1$ such that $k(2
\g-1)<1+\epsilon/5$.  From now on, once a given $\g$ close to $1$ is chosen,
$\epsilon$ will always denote a small number with this property.

\subsection{Capacities and trees}

The $\g$-capacity of a set $X$ in an interval $I$ is defined as
follows:
$$
p_\g(X|I)=\sup_{h \in QS(\g)} \frac {|h(X \cap I)|} {|h(I)|}.
$$

This geometric quantity is well adapted to our context, since it is well
behaved under tree decompositions of sets.  In other words, if $I^j$ are
disjoint subintervals of $I$ and $X \subset \cup I^j$ then
$$
p_\g(X|I) \leq p_\g(\cup_j I_j|I) \sup_j p_\g(X|I^j).
$$

\subsection{A measure-theoretical lemma}

Our procedure consists in obtaining successively smaller (but still
full-measure) classes of maps for which we can give a progressively
refined statistical description of the dynamics.  This is done inductively
as follows: we pick a class $X$ of maps (which we have previously shown to
have full measure among non-regular maps) and for each map in $X$
we proceed to describe the dynamics (focusing on the statistical behavior
of return and landing maps for deep levels of the principal nest),
then we use this information to show that a subset $Y$ of $X$
(corresponding to parameters for which the statistical behavior of the
{\it critical orbit} is not anomalous) still has full measure.
An example of this parameter exclusion process is done by Lyubich in \cite
{parapuzzle} where he shows using a probabilistic argument that
the class of simple maps has full measure in $\FF$.

Let us now describe our usual argument (based on the argument of Lyubich
which in turn is a variation of the Borel-Cantelli Lemma).
Assume at some point we know how to prove that almost every simple map
belongs to a certain set $X$.
Let $Q_n$ be a (bad) property that a map
may have (usually some anomalous statistical parameter
related to the $n$-th stage of the principle nest).
Suppose we prove that if $f \in X$ then
the probability that a map in $J_n(f)$ has the property $Q_n$ is bounded by
$q_n(f)$ which is shown to be summable
for all $f \in X$.  We then conclude that almost every map does not have
property $Q_n$ for $n$ big enough.

Sometimes we also apply the same argument, proving instead that $q_n(f)$
is summable where $q_n(f)$ is the
probability that a map in $J^{\tau_n}_n(f)$ has property $Q_n$,
(recall that $\tau_n$ is such that $f \in J^{\tau_n}_n(f)$).

In other words, we apply the following general result.

\begin{lemma} \label {measureth}

Let $X \subset \R$ be a measurable set such that for each
$x \in X$ is defined a sequence
$D_n(x)$ of nested intervals converging to $x$
such that for all $x_1,x_2 \in X$
and any $n$, $D_n(x_1)$ is either equal or disjoint to $D_n(x_2)$.  Let
$Q_n$ be measurable subsets of
$\R$ and $q_n(x)=|Q_n \cap D_n(x)|/|D_n(x)|$.  Let $Y$ be
the set of all $x \in X$ which belong to at most
finitely many $Q_n$.
If $\sum q_n(x)$ is finite for almost any $x \in X$ then $|Y|=|X|$.

\end{lemma}

\begin{pf}

Let $Y_n=\{x \in X|\sum_{k=n}^\infty q_k(x)<1/2\}$.  It is clear that $Y_n
\subset Y_{n+1}$ and $|\cup Y_n|=|X|$.

Let $Z_n=\{x \in Y_n||Y_n \cap D_m(x)|/|D_m(x)|>1/2, m \geq n\}$.
It is clear that $Z_n \subset Z_{n+1}$ and $|\cup Z_n|=|X|$.

For $m \geq n$,
let $T^m_n=\cup_{x \in Z_n} D_m(x)$.  Let $K^m_n=T^m_n \cap Q_m$.
Of course
$$
|K^m_n|=\int_{T^m_n} q_m \leq 2 \int_{Y_n} q_m.
$$

And of course
$$
\sum_{m \geq n} \int_{Y_n} q_m \leq \frac{1}{2} |Y_n|.
$$

This shows that $\sum_{m \geq n} |K^m_n| \leq |Y_n|$,
so almost every point in $Z_n$ belongs to at most finitely many $K^m_n$.
We conclude then that almost every point in
$X$ belongs to at most finitely many $Q_m$.
\end{pf}

The following obvious reformulation will be often convenient

\begin{lemma} \label {simp}

In the same context as above, assume that we are given sequences $Q_{n,m}$,
$m \geq n$ of measurable sets and let $Y_n$ be the set of $x$ belonging to
at most finitely many $Q_{n,m}$.  Let $q_{n,m}(x)=|Q_{n,m} \cap
D_m(x)|/|D_m(x)|$.
Let $n_0(x) \in \N \cup \{\infty\}$ be such that
$\sum_{m=n}^\infty q_{n,m}(x)<\infty$ for $n \geq n_0(x)$. 
Then for almost every $x \in X$, $x \in Y_n$ for $n \geq n_0(x)$.

\end{lemma}

In practice, we will estimate the capacity of sets in the phase space:
that is, given a map $f$ we will obtain subsets $\tilde Q_n[f]$ in the phase
space, corresponding to bad branches of return or landing maps.  We will
then show that for some $\g>1$ we have
$\sum p_\g(\tilde Q_n[f]|I_n[f])<\infty$ or
$\sum p_\g(\tilde Q_n[f]|I^{\tau_n}_n[f])<\infty$.  We will then use PhPa2
or PhPa1, and the measure-theoretical lemma above to conclude that with
total probability among non-regular maps, for all $n$ sufficiently big,
$R_n(0)$ does not belong to a bad set.

From now on when we prove that almost every non-regular map
has some property, we will just say that with total probability
(without specifying) such property holds.

(To be strictly formal, we have fixed the renormalization level
$\kappa$ (in particular to define the sequence $J_n$ without ambiguity),
so applications of the measure theoretical argument will actually
be used to conclude that for almost every parameter in
$\FF_\kappa$ a given property holds.  Since almost every
non-regular map belongs
to some $\FF_k$, this is equivalent to the statement regarding almost every
non-regular parameter.)

\section{Statistics of the principal nest}

\subsection{Decay of geometry}

Let as before $\tau_n \in \Z$ be such that $R_n(0) \in I^{\tau_n}_n$.

An important parameter in our construction will be the scaling factor
$$
c_n=\frac {|I_{n+1}|} {|I_n|}.
$$
This variable of course changes inside each
$J^{\tau_n}_n$ window, however, not by much.
From PhPh1, for instance, we get
that with total probability
$$
\lim_{n \to \infty} \sup_{g_1,g_2 \in J^{\tau_n}_n} \frac {\ln(c_n[g_1])}
{\ln(c_n[g_2])}=1.
$$

This variable is by far the most important on our analysis of the statistics
of return maps.  We will often consider other variables
(say, return times): we will show that the distribution of those
variables is concentrated near some average value.
Our estimates will usually give a range
of values near the average, and $c_n$ will play an important role.  Due
(among other issues) to the variability of $c_n$ inside the parameter
windows, the ranges we select will depend on $c_n$ up to an exponent
(say, between $1-\epsilon$ and $1+\epsilon$), where $\epsilon$ is a small, but fixed,
number.  From the estimate we just obtained,
for big $n$ the variability (margin of error) of $c_n$ will fall
comfortably in such range, and we won't elaborate more.

A general estimate on the rates decay of $c_n$ was obtained by Lyubich: he
shows that (for a finitely
renormalizable unimodal map with a recurrent critical point), $c_{n_k}$
decays exponentially (on $k$), where $n_k-1$ is the subsequence of
non-central levels of $f$.  For simple maps, the same is true with
$n_k=k$, as there are only finitely many central returns.
Thus we can state:

\begin{thm}[see \cite {attractors}] \label {attractors theorem}

If $f$ is a simple map then there exists $C>0$, $\lambda<1$ such that $c_n<C
\lambda^n$.

\end{thm}

Let us use the following notation for the combinatorics of a point $x \in
I_n$.  If $x \in I^j_n$ we let $j^{(n)}(x)=j$ and if
$x \in C^\d_n$ we let $\d^{(n)}(x)=\d$.

\begin{lemma} \label {estimate on m}
  
With total probability, for all $n$ sufficiently big we have
\begin{align}
\label {est1}
&p_{2 \g-1}(|\d^{(n)}(x)| \leq k|I_n) <
k c_n^{1-\epsilon/2},\\
\label {est2}
&p_{2 \g-1}(|\d^{(n)}(x)| \geq k|I_n) <
e^{-k c_n^{1+\epsilon/2}}.
\end{align}  
We also have
\begin{align}
&p_{2 \g-1}(|\d^{(n)}(x)| \leq k|I^{\tau_n}_n) <
k c_n^{1-\epsilon/2},\\
&p_{2 \g-1}(|\d^{(n)}(x)| \geq k|I^{\tau_n}_n) <
e^{-k c_n^{1+\epsilon/2}}.
\end{align}
  
\end{lemma}

\begin{pf}

Let us compute the first two estimates.

Since $I^0_n$ is in the middle of
$I_n$, we have as a simple consequence of the Real Schwarz Lemma (see \cite
{attractors} and (\ref {8.3}) in Lemma \ref {rsl} below) that
$$
\frac {c_n}{4}<\frac {|C^\d_n|}{|I^\d_n|}<4 c_n.
$$
As a consequence
$$
p_{2 \g-1}(|\d^{(n)}(x)|=m|I_n)<(4 c_n)^{1-\epsilon/3}
$$
and we get the estimate (\ref {est1}) summing up on $0 \leq m \leq k$.

For the same reason, we get that
\begin{align*}
p_{2 \g-1}(|\d^{(n)}&(x)| \geq m+1|I_n)\\
&<\left (1-\left (\frac {c_n} {4}\right )^{1+\epsilon/3} \right)
p_{2 \g-1}(|\d^{(n)}(x)| \geq m|I_n).
\end{align*}
This implies
$$
p_{2 \g-1}(|\d^{(n)}(x)| \geq m|I_n)
\leq \left (1-\left (\frac {c_n} {4}
\right )^{1+\epsilon/3}\right )^m .
$$
  
Estimate (\ref {est2}) follows from
\begin{align*}
\left (1-\left (\frac {c_n} {4}\right )^{1+\epsilon/3}\right )^k &<
(1-c_n^{1+\epsilon/2})^k\\
&<((1-c_n^{1+\epsilon/2})^{c_n^{-1-\epsilon/2}})^{k
c_n^{1+\epsilon/2}}\\
&<e^{-k c_n^{1+\epsilon/2}}.
\end{align*} 
   
The two remaining estimates are analogous.
\end{pf}

Let us now transfer this result
(more precisely the second pair of estimates)
to the parameter in each
$J^{\tau_n}_n$ window using PhPa1.  To do this notice
that the measure of the
complement of the set of parameters in $J^{\tau_n}_n$ such that
$c_n^{-1+2\epsilon}<s_n<c_n^{-1-2\epsilon}$ can be estimated by
$2 c_n^\epsilon$ for $n$ big which is summable
for all $\epsilon$ by Theorem \ref {attractors theorem}.  So we have:

\begin{lemma} \label {growth of s_n}

With total probability,
$$
\lim_{n \to \infty} \frac {\ln (s_n)} {\ln (c^{-1}_n)}=1.
$$
\end{lemma}

The parameter $s_n$ influences the size of $c_{n+1}$ in a determinant way.

\begin{cor} \label {c_n torrential}

With total probability,
\be \label {c_n tor}
\liminf_{n \to \infty}
\frac {\ln(\ln(c_{n+1}^{-1}))} {\ln(c_n^{-1})} \geq 1.
\ee
In particular, $c_n$
decreases at least torrentially fast.

\end{cor}

\begin{pf}

It is easy to see (using for instance the Real Schwarz Lemma,
see \cite {attractors}, see also item (\ref {8.4}) in
Lemma \ref {rsl} below)
that there exists a constant $K>0$ (independent of $n$) such that
for each $\d \in \Omega$, both components of
$I^{\sigma^+(\d)}_n \setminus I^\d_n$ have size at least
$(e^K-1)|I^\d_n|$.  In particular, by induction,
if $R_n(0) \in C^\d_n$ we have that both
gaps of $I_n \setminus C^\d_n$ have size at least $(e^{K s_n}-1) |C^\d_n|$. 
Taking the preimage by $R_n$, and using the Real Schwarz Lemma again, we see
that $c_{n+1}<C e^{K s_n/2}$ for some constant $C>0$ independent of $n$.  We
conclude that
$$
\liminf \frac {\ln (c_{n+1}^{-1})} {s_n} \geq \frac {K} {2},
$$
and since $c_n \to 0$ as $n \to \infty$ we have
$$
\liminf \frac {\ln (\ln (c_{n+1}^{-1}))} {\ln (s_n)} \geq 1
$$
which together with Lemma \ref {growth of s_n} implies (\ref {c_n tor}).
\end{pf}

\begin{rem}

In the proof of Corollary \ref {c_n torrential}, the constant $K>0$
is related to the real bounds.  In our situation,
since we have decay of geometry, we can actually take $K \to \infty$ as
$n \to \infty$, so we actually have
$$
\frac {\ln (c_{n+1}^{-1})} {s_n} \to \infty
$$
torrentially fast.

\end{rem}

\subsection{Fine partitions}

We use Cantor sets $K_n$ and $\tilde K_n$ to partition the phase space.
In many circumstances we are directly concerned with intervals of this
partition.  However, sometimes we just want to exclude an interval of given
size (usually a neighborhood of $0$).
This size does not usually correspond to (the closure of) a union
of gaps, so we instead should consider in applications an interval which
is a union of gaps, with approximately the given size
\footnote{
We need to consider intervals which are union of gaps due to our phrasing of
the Phase-Parameter relation, which only gives information about such gaps. 
However, this is not absolutely necessary, and we could have proceeded in a
different way: the proof of the Phase-Parameter relation actually shows that
there is a {\it holonomy map} between phase and parameter
{\it intervals} (and not only Cantor sets) corresponding to a {\it
holomorphic motion} for which we can obtain good qs estimates.  While this
map is not canonical, the fact that it is a holonomy map for a
holomorphic motion with good qs estimates would allow our
proofs to work.}.
The degree of relative approximation will always be torrentially
good (in $n$), so we
usually won't elaborate on this.  In this section we just give some
results which will
imply that the partition induced by the Cantor sets are fine enough to allow
torrentially good approximations.

The following lemma summarizes the situation.  The proof is based on
estimates of distortion using the Real Schwarz Lemma and
the Koebe Principle (see \cite {attractors}) and is very simple,
so we just sketch the proof.

\begin{lemma} \label {rsl}

The following estimates hold:

\begin{align}
\label {8.1}
\frac {|I^j_n|} {|I_n|}=O(\sqrt {c_{n-1}}),\\
\label {8.2}
\frac {|I^\d_n|} {|I^{\sigma^+(\d)}_n|}=O(\sqrt {c_{n-1}}),\\
\label {8.3}
\frac {c_n} {4}<\frac {|C^{\d}_n|} {|I^{\d}_n|}<4 c_n,\\
\label {8.4}
\frac {|\tilde I_{n+1}|} {|I_n|}=O(e^{-s_{n-1}}).
\end{align}

\end{lemma}

\begin{pf}

(Sketch.)  Since $R^\d_n$ has negative Schwarzian derivative, it immediately
follows that the Koebe space\footnote {The Koebe space of an interval $T'$
inside an interval $T \supset T'$ is the minimum of $|L|/|T'|$ and
$|R|/|T'|$ where $L$ and $R$
are the components of $T \setminus T'$.  If the Koebe
space of $T'$ inside $T$ is big, then the Koebe
Principle states that a diffeomorphism onto $T'$ which has an
extension with negative Schwarzian derivative onto
$T$ has small distortion.  In this case, it follows that the
Koebe space of the preimage of $T'$ inside the preimage of $T$
is also big.} of $C^\d_n$ inside $I^\d_n$ has at least order $c_n^{-1}$.

It is easy to see that $R_{n-1}|_{I_n}$ can
be written as $\phi \circ f$ where $\phi$ extends to a diffeomorphism onto
$I_{n-2}$ with negative
Schwarzian derivative and thus with very small distortion.  Since
$R_{n-1}(I^j_n)$ is contained on some $C^\d_{n-1}$, we see that the Koebe
space of $I^j_n$ in $I_n$ is at least of order $c_{n-1}^{-1/2}$ which
implies (\ref {8.1}).

Let us now consider an interval $I^\d_n$.  Let $I^j_n$ be such that
$R_n^{\sigma^+(\d)}(I^\d_n)=I^j_n$.  We can pullback the Koebe space
of $I^j_n$ inside $I_n$ by $R_n^{\sigma^+(\d)}$, so (\ref {8.1}) implies
(\ref {8.2}).  Moreover, this shows by induction that
the Koebe space of $I^\d_n$ inside
$I_n$ is at least of order $c_{n-1}^{-|\d|/2}$.
Since $R_{n-1}(\tilde I_{n+1}) \subset
I^\d_{n-1}$ with $|\d|=s_{n-1}$,
the Koebe space of $\tilde I_{n+1}$ in $I_n$ is at
least $c_{n-2}^{-|\d|/4}$, which implies (\ref {8.4}).

It is easy to see
that $R^\d_n|_{I^\d_n}$ can be written as $\phi \circ f \circ
R^{\sigma^+(\d)}_n$, where $\phi$ has small distortion.  Due to (\ref
{8.1}), $R^{\sigma^+(\d)}_n|_{I^\d_n}$ also has small distortion, so a
direct computation with $f$ (which is purely quadratic) gives (\ref {8.3}).
\end{pf}

In other words, distances in $I_n$ can be measured with precision $\sqrt
{c_{n-1}}|I_n|$ in the partition induced by $\tilde K_n$, due to
(\ref {8.1}) and (\ref {8.4}) (since $e^{-s_{n-1}} \ll c_{n-1}$).

Distances can be measured much more precisely with respect to the partition
induced by $K_n$, in fact we have good precision in each
$I^\d_n$ scale.  In other words, inside $I^\d_n$,
the central gap $C^\d_n$ is
of size $O(c_n|I^\d_n|)$ (by (\ref {8.3}))
and the other gaps have size $O(\sqrt {c_{n-1}} |C^{\d}_n|)$
(by (\ref {8.2}) and (\ref {8.3})).

\subsection{Initial estimates on distortion}

To deal with the distortion control we need some preliminary known results.
Those estimates are based on the Koebe Principle and the
estimates of Lemma \ref {rsl}.  All needed arguments are already
contained in the proof of Lemma \ref {rsl}, so we won't get into details.

\begin{prop}

The following estimates hold:

\begin{enumerate}

\item For any $j$, if $R_n|_{I^j_n}=f^k$,
$\dist(f^{k-1}|_{f(I^j_n)})=1+O(c_{n-1})$,

\item For any $\d$,
$\dist(R^{\sigma^+(\d)}_n|_{I^\d_n})=1+O(\sqrt {c_{n-1}})$.

\end{enumerate}

\end{prop}

We will use the following immediate consequence for the decomposition of
certain branches.

\begin{lemma} \label {decomposition}

With total probability,

\begin{enumerate}

\item $R_n|_{I^0_n}=\phi \circ f$ where $\phi$ has torrentially small
distortion,

\item $R^{\d}_n=\phi_2 \circ f \circ \phi_1$ where $\phi_2$ and
$\phi_1$ have torrentially small distortion and
$\phi_1=R^{\sigma^+(\d)}_n$.

\end{enumerate}

\end{lemma}

\subsection{Estimating derivatives}

\begin{lemma} \label {away from the boundary}

Let $w_n$ denote the relative distance in
$I_n$ of $R_n(0)$ to $\partial I_n \cup \{0\}$:
$$
w_n=\frac {d(R_n(0),\partial I_n \cup \{0\})} {|I_n|}, \quad \text {where }
d(x,X)=\inf_{y \in X} |y-x|.
$$
With total probability,
$$
\limsup_{n \to \infty} \frac {-\ln(w_n)} {\ln(n)} \leq 1.
$$
In particular $R_n(0) \notin \tilde I_{n+1}$ for all $n$ large enough.

\end{lemma}

\begin{pf}

This is a simple consequence of PhPa2, using that
$n^{-1-\delta}$ is summable, for all $\delta>0$ (using (\ref {8.4}) to
obtain the last conclusion).
\end{pf}

From now on we suppose that $f$ satisfies the conclusions of the above
lemma.

\begin{lemma} \label {dist}

With total probability,
$$
\limsup_{n \to \infty}
\frac {\sup_{j \neq 0} \ln(\dist(f|_{I^j_n}))} {\ln(n)} \leq 1/2.
$$

\end{lemma}

\begin{pf}

Denote by $P^{\d}_n$ a $|C^{\d}_n|/n^{1+\delta}$
neighborhood of $C^{\d}_n$.  Notice that the gaps of the Cantor sets
$K_n$ inside $I^\d_n$ which are different from $C^\d_n$
are torrentially (in $n$) smaller then $C^\d_n$, so we can take
$P^\d_n$ as a union of gaps of $K_n$ up to torrentially small error.

It is clear that if $h$ is a $\g$-qs homeomorphism ($\g$ close to $1$)
then
$$
|h(P^\d_n \setminus C^\d_n)| \leq n^{-1-\delta/2}|h(C^\d_n)|
$$
Notice that if $C^\d_n$ is contained in $I^j_n$ with $j \neq \tau_n$,
then $P^\d_n$ does not intersect $I^{\tau_n}_n$.
Since the $C^\d_n$ are disjoint,
$$
p_\g(\cup (P^\d_n \setminus C^\d_n)|I^{\tau_n}_n) \leq
n^{-1-\delta/2}
$$
which is summable.

Transferring this estimate to the parameter using
PhPa1 we see that with total probability, if $n$ is sufficiently big,
if $R_n(0)$ does not belong to
$C^\d_n$ then $R_n(0)$ does not belong to $P^\d_n$ as well.
In particular, if $n$ is sufficiently big, the critical point $0$ will
never be in a $n^{-1/2-\delta/5}|I^j_{n+1}|$ neighborhood of any $I^j_{n+1}$ with
$j \neq 0$ (the change from $n^{-1-\delta}$ to $n^{-1/2-\delta/5}$
is due to taking the inverse image by $R_n|_{I_{n+1}}$, which
corresponds, up to torrentially small distortion,
to taking a square root, and causes the
division of the exponent by two).  This implies the required estimate on
distortion since $f$ is quadratic.
\end{pf}

\begin{lemma} \label {distortion}

With total probability,
\begin{equation} \label {disto}
\limsup_{n \to \infty} \frac {\sup_{\d \in \Omega} \ln(\dist(R^\d_n))}
{\ln(n)} \leq \frac {1} {2}.
\end{equation}
In particular, for $n$ big enough,
$\sup_{\d \in \Omega} \dist(R^\d_n) \leq 2^n$ and
$|DR_n(x)|>2$, $x \in \cup_{j \neq 0} I^j_n$.

\end{lemma}

\begin{pf}

By Lemma \ref {decomposition}, Lemma \ref {dist} implies (\ref {disto}).  If
$j \neq 0$, by (\ref {8.1}) of Lemma \ref {rsl} we get that
$|R_n(I^j_n)|/|I^j_n|=|I_n|/|I^j_n|>c_{n-1}^{-1/3}$, so
$\dist(R_n|_{I^j_n}) \leq 2^n$ implies that for all $x \in I^j_n$,
$|DR_n(x)|>c_{n-1}^{-1/3} 2^{-n}>2$.
\end{pf}

\begin{rem} \label {approx}

Lemma \ref {dist} has also an application for approximation of
intervals.  This result implies that if $I^j_n=(a,b)$ and $j \neq 0$,
we have $1/2^n<b/a<2^n$.  As a consequence, for any symmetric (about $0$)
interval $I_{n+1} \subset X \subset I_n$, there exists a symmetric
(about $0$) interval $X \subset \tilde X$, which is union of
$I^j_n$ and such that $|\tilde X|/|X|<2^n$ (approximation by union of
$C^\d_n$, with $|\tilde X|/|X|$ torrentially close to $1$, follows
more easily from the discussion on fine partitions).

\end{rem}

We will also need to estimate derivatives of iterates of $f$, and not only
of return branches.

\begin{lemma} \label {lower bound}

With total probability, if $n$ is sufficiently big and if
$x \in I^j_n$, $j \neq 0$, and $R_n|_{I^j_n}=f^K$, then for $1 \leq k
\leq K$, $|(Df^k(x))|>|x| c_{n-1}^3$.

\end{lemma}

\begin{pf}

First notice that by Lemma \ref {away from the boundary} and Lemma
\ref {decomposition},
$R_n|_{I^0_n}=\phi \circ f$ with $|D\phi|>1$, provided $n$ is big enough
(since $\phi$ has small distortion and there is a big macroscopic expansion
from $f(I^0_n)$ to $R_n(I^0_n)$).
Also, by Lemma \ref {c_n torrential}, $|I_n|$ decays so fast that
$\prod_{r=1}^n |I_n|>c_{n-1}^{3/2}$
for $n$ big enough.
Finally, by Lemma \ref {distortion}, for $n$ big enough, $|DR_n(x)|>1$ for
$x \in I^j_n$, $j \neq 0$.
Let $n_0$ be so big that if $n \geq n_0$, all the above properties hold.

From hyperbolicity of $f$
restricted to the complement of $I_{n_0}$ (from Lemma \ref {hyperbol}),
there exists a constant $C>0$ such that if $s_0$ is such that
$f^s(x) \notin I^0_{n_0}$ for every $s_0 \leq s <k$ then
$|Df^{k-s_0}(f^{s_0}(x))|>C$.

Let us now consider some $n \geq n_0$.
If $k=K$, we have a full return and
the result follows from Lemma \ref {distortion}.

Assume now $k<K$.  Let us define $d(s)$,
$0 \leq s \leq k$ such that $f^s(x) \in I_{d(s)}
\setminus I^0_{d(s)}$ (if $f^s(x) \notin I_0$ we set $d(s)=-1$).
Let $m(s)=\max_{s \leq t \leq k} d(t)$.  Let us define a finite sequence
$\{k_r\}_{r=0}^l$ as follows.  We let $k_0=0$ and
supposing $k_r<k$ we let $k_{r+1}=\max \{k_r<s \leq k|
d(s)=m(s)\}$.
Notice that $d(k_i)<n$ if $i \geq 1$, since otherwise
$f^{k_i}(x) \in I_n$ so $k=k_i=K$ which contradicts our
assumption.

The sequence $0=k_0 < k_1<...<k_l=k$ satisfies
$n=d(k_0)>d(k_1)>...>d(k_l)$.
Let $\theta$ be maximal with $d(k_\theta) \geq n_0$.  We have of course
$$
|Df^{k-k_\theta}(f^{k_\theta}(x))|>C|Df(f^{k_\theta}(x))|,
$$
so if $\theta=0$ then $Df^k(x)>|2 C x|$ and we are done.

Assume now $\theta>0$.  We have of course
$$
|Df^{k-k_\theta}(f^{k_\theta}(x))| >
C|Df(f^{k_\theta}(x))|>C|I_{d(k_\theta)+1}|
$$

For $1 \leq r \leq \theta$, the action of $f^{k_r-k_{r-1}}$ near
$f^{k_{r-1}}(x)$ is obtained by applying
the central component of $R_{d(k_r)}$ followed by several non-central
components of $R_{d(k_r)}$.  Since $d(k_r) \geq n_0$, we can estimate
$$
|Df^{k_r-k_{r-1}}(f^{k_{r-1}}(x))|>|DR_{d(k_r)}(f^{k_{r-1}}(x))| >
|Df(f^{k_{r-1}}(x))|.
$$
For $r=1$, this argument gives
$|Df^{k_1}(x)| \geq |Df(x)|$, while for $r>1$ we can estimate
$$
|Df^{k_r-k_{r-1}}(f^{k_{r-1}}(x))|>|Df(f^{k_{r-1}}(x))| >
|I_{d(k_{r-1})+1}|.
$$

Combining it all we get
\begin{align*}
|Df^k(x)|&=|Df^{k_1}(x)| \cdot |Df^{k-k_\theta}(f^{k_\theta}(x))|
\prod_{r=2}^\theta |Df^{k_r-k_{r-1}}(f^{k_{r-1}}(x))|\\
&>|2 x| \cdot C \cdot |I_{d(k_\theta)+1}|
\prod_{r=2}^\theta |I_{d(k_{r-1})+1}| =
|2 C x| \prod_{r=1}^\theta |I_{d(k_r)+1}|\\
&\geq |2 C x| \prod_{r=0}^n |I_r|>|x| c_{n-1}^3.
\end{align*}
\end{pf}

\section{Sequences of quasisymmetric constants and trees} \label {seq
constants}

\subsection{Preliminary estimates}

From now on, we will need to transfer estimates on the capacity of certain
sets from level to level of the principal nest.  In order to do so
we will need to consider not only $\g$-capacities with some
$\gamma$ fixed, but different constants for different levels of the
principal nest.  To do so, we will make use of sequences of constants
converging (decreasing) to a given value $\g$.  We recall that $\g$ is some
constant very close to $1$ such that
$k(2 \gamma-1)<1+\epsilon/5$, with $\epsilon$ very small.

We define the sequences $\r_n=(n+1)/n$ and $\tr_n=(2n+3)/(2n+1)$,
so that $\r_n>\tr_n>\r_{n+1}$ and $\lim \r_n=1$.
We define the sequence $\g_n=\g \r_n$ and an intermediate sequence
$\tg_n=\g \tr_n$.

As we know, the generalized renormalization process relating $R_n$ to
$R_{n+1}$ has two phases, first
$R_n$ to $L_n$ and then $L_n$ to $R_{n+1}$.  The following
remarks shows why it is useful to consider the sequence
of quasisymmetric constants due to losses related to distortion.

\begin{rem} \label {remark 1}

Let $S$ be an interval contained in $I^\d_n$.
Using Lemma \ref {decomposition}
we have $R^\d_n|_S=\psi_2 \circ f \circ \psi_1$,
where the distortion of $\psi_2$ and $\psi_1$ are torrentially small
and $\psi_1(S)$ is contained in some $I^j_n$, $j \neq 0$.
If $S$ is contained in $I^0_n$ we may as well write $R_n|_S=\phi \circ f$,
and the distortion of $\phi$ is also torrentially small.

In either case, if we decompose $S$ in $2km$ intervals
$S_i$ of equal length, where $k$ is the distortion of either
$R^\d_n|_S$ or $R_n|_S$ and $m$ is subtorrentially big
(say, $m<2^n$), the distortion obtained restricting to any interval $S_i$
will be bounded by $1+m^{-1}$.  Indeed, in the case $S \subset I^0_n$,
we have $\dist (R_n|_{S_i}) \leq \dist(\phi) \dist(f|_{S_i})$.
Now $k=\dist(R_n|_S) \geq \dist(\phi)^{-1} \dist(f|_S)$.  Since $f$ is
quadratic,
\begin{align}
\nonumber
\dist(f|_{S_i})-1 &\leq \frac {|S_i|} {|S|} (\dist(f|_S)-1) \leq
\frac {1} {2km} (k \dist(\phi)-1)\\
\nonumber
&\leq \frac {\dist (\phi)} {2m}.
\end{align}
Since $\dist(\phi)-1$ is torrentially small, $\dist(f|_{S_i}) \leq
1+(2/3)m^{-1}$ and $\dist(R_n|_{S_i}) \leq 1+m^{-1}$.
The case $S \subset I^\d_n$ is entirely
analogous, considering $\dist(R^\d_n|_{S_i}) \leq
\dist(\psi_2)\dist(f|_{\psi_1(S_i)})\dist(\psi_1)$, and
using torrentially small distortion of $\psi_1$ and $\psi_2$.
The estimate now becomes
\begin{align}
\nonumber
\dist(f|_{\psi_1(S_i)})-1 &\leq \frac {|\psi_1(S_i)|} {|\psi_1(S)|}
(\dist(f|_{\psi_1(S)})-1) \leq
\frac {\dist(\psi_1)} {2km} (k \dist(\psi_1)\dist(\psi_2)-1)\\
\nonumber
&\leq \frac {\dist (\psi_1)^2 \dist(\psi_2))} {2m}
\end{align}
and we conclude again that $\dist(R^\d_n|_{S_i}) \leq 1+m^{-1}$.

\end{rem}

\begin{rem} \label {remark 2}

Let us fix now $\g$ such that the corresponding $\epsilon$ is small enough. 
We have the following estimate for the effect of the pullback of a subset of
$I_n$ by the central branch $R_n|_{I^0_n}$.
With total probability, for all $n$ sufficiently big,
if $X \subset I_n$ satisfies
$$
p_{\tg_n}(X|I_n)<\delta \leq n^{-1000}
$$
then
$$
p_{\g_{n+1}}((R_n|_{I_{n+1}})^{-1}(X)|I_{n+1})<\delta^{1/5}.
$$

Indeed, let $V$ be a $\delta^{1/4}|I_{n+1}|$ neighborhood of $0$.  Then
$R_n|_{I_{n+1} \setminus V}$ has distortion bounded by $2 \delta^{1/4}$.

Let $W \subset I_n$ be an interval of size $\lambda |I_n|$.  Of course
$$
p_{\tg_n}(X \cap W|W)<\delta \lambda^{-1-\epsilon}.
$$

Let us decompose each side of $I_{n+1} \setminus V$ as a union of
$n^3 \delta^{-1/4}$ intervals of equal length.
Let $W$ be such an interval.
From Lemma \ref {away from the boundary},
it is clear that the image of $W$
covers at least $\delta^{1/2} n^{-4} |I_n|$.
It is clear then that
$$
p_{\tg_n}(X \cap R_n(W)|R_n(W))<\delta^{(1-\epsilon)/2} n^{4+4\epsilon}.
$$
So we conclude that (since the distortion of $R_n|_W$ is of order
$1+n^{-3}$ by Remark \ref {remark 1})
$$
p_{\g_{n+1}}((R_n|_{I_{n+1}})^{-1}(X) \cap W|W)<
\delta^{(1-\epsilon)/2} n^5
$$
(we use the fact that the composition of a $\g_{n+1}$-qs map with a map with
small distortion in $\tg_n$-qs).
Since
$$
p_{\g_{n+1}}(V|I_{n+1})<(2 \delta^{1/4})^{1-\epsilon},
$$
we get the required estimate.

\end{rem}

\subsection{More on trees} \label {more tree}

Let us see an application of the above remarks.

\begin{lemma} \label {2n}

With total probability, for all $n$ sufficiently big
$$
p_{\tg_n}((R^\d_n)^{-1}(X)|I^\d_n)<2^n p_{\g_n}(X|I_n).
$$

\end{lemma}

\begin{pf}

Decompose $I^\d_n$ in $n^{\ln(n)}$ intervals of equal length, say,
$\{W_i\}_{i=1}^{n^{\ln(n)}}$.  Then by Lemma \ref {distortion},
$|R_n^\d(W_i)|>n^{-2 \ln n} |I_n|$, so we get
$$
p_{\g_n}(R^\d_n(W_i) \cap X|R^\d_n(W_i))<n^{4 \ln(n)} p_{\g_n}(X|I_n).
$$

Applying Remark \ref {remark 1}, we see that 
$$
p_{\tg_n}((R^\d_n)^{-1}(X) \cap W_i|W_i)<n^{4 \ln (n)} p_{\g_n}(X|I_n),
$$
(we use the fact that the composition of a $\tg_n$-qs map with a map with
small distortion is $\g_n$-qs) which implies the desired estimate.
\end{pf}

By induction we get:

\begin{lemma} \label {distortion product}

With total probability, for $n$ is big enough, if
$X_1,...,X_m \subset \Z \setminus \{0\}$
\begin{align*}
p_{\tg_n}(\d^{(n)}(x)=(j_1,...,j_m,...,j_{|\d^{(n)}(x)|}),&
j_i \in X_i, 1 \leq i \leq m|I_n)\\
&\leq 2^{mn} \prod_{i=1}^m p_{\g_n}(j^{(n)}(x) \in X_i|I_n).
\end{align*}

\end{lemma}

The following is an obvious variation of the previous lemma
fixing the start of the sequence.

\begin{lemma} \label {distortion product2}

With total probability, for $n$ is big enough, if
$X_1,...,X_m \subset \Z \setminus \{0\}$, and if $\d=(j_1,...,j_k)$
we have
\begin{align*}
p_{\tg_n}&(\d^{(n)}(x)=
(j_1,...,j_k,j_{k+1},...,j_{k+m},...,j_{|\d^{(n)}(x)|}),
j_{i+k} \in X_i, 1 \leq i \leq m|I^\d_n)\\
&\leq 2^{mn} \prod_{i=1}^m p_{\g_n}(j^{(n)}(x) \in X_i|I_n).
\end{align*}
In particular, with $\d=(\tau_n)$,
\begin{align*}
p_{\tg_n}&(\d^{(n)}(x)=(\tau_n,j_1,...,j_m,j_{m+1},...,j_{|\d^{(n)}(x)|}),
j_i \in X_i, 1 \leq i \leq m|I^{\tau_n}_n)\\
&\leq 2^{mn} \prod_{i=1}^m p_{\g_n}(j^{(n)}(x) \in X_i|I_n).
\end{align*}

\end{lemma}

The last part of the above lemma will be often necessary
in order to apply PhPa1.

Sometimes we are more interested in the case where the $X_i$ are all equal.

Let $Q \subset \Z \setminus \{0\}$.  Let $Q(m,k)$ denote the set of
$\d=(j_1,...,j_m)$ such that $\#\{1 \leq i \leq m,\, j_i \in Q\} \geq k$.

Define $q_n(m,k)=
p_{\tg_n}(\cup_{\d \in Q(m,k)} I^{\d}_n|I_n)$.

Let $q_n=p_{\g_n}(\cup_{j \in Q} I^j_n|I_n)$.

\begin{lemma} \label {tree estimate}

With total probability, for $n$ large enough,
\begin{equation} \label {q_n}
q_n(m,k) \leq \binom {m} {k} (2^n q_n)^k.
\end{equation}

\end{lemma}

\begin{pf}

We have the following recursive estimates for $q_n(m,k)$:
\begin{enumerate}

\item $q_n(1,0)=1$, $q_n(1,1) \leq q_n \leq 2^n q_n$, and
$q_n(m+1,0) \leq 1$ for $m \geq 1$,

\item $q_n(m+1,k+1) \leq q_n(m,k+1)+2^n q_n q_n(m,k)$.

\end{enumerate}

Indeed, (1) is completely obvious and
if $(j_1,...,j_{m+1}) \in Q(m+1,k+1)$ then either $(j_1,...,j_m) \in
Q(m,k+1)$ or $(j_1,...,j_m) \in Q(m,k)$ and $j_{m+1} \in Q$, so (2)
follows from Lemma \ref {2n}.  It is clear that (1) and (2) imply by
induction (\ref {q_n}).
\end{pf}

We recall that by Stirling Formula,
$$
\binom {m} {q m}<\frac {m^{q m}} {(q m)!} <
\left (\frac {3} {q} \right )^{q m}.
$$

So we can get the following estimate.  For $q \geq q_n$,
\begin{equation} \label {q}
q_n(m,(6 \cdot 2^n) q m)<\left (\frac {1} {2}\right )^{(6 \cdot 2^n) q m}.
\end{equation}

It is also used in the following form.  If $q^{-1}>6 \cdot 2^n$ (it is
usually the case, since $q$ will be torrentially small)
\begin{equation} \label {sumq}
\sum_{k>q^{-2}} q_n(k,(6 \cdot 2^n) q k) <
2^{-n} q^{-1} \left (\frac {1} {2} \right )^{(6 \cdot 2^n) q^{-1}}.
\end{equation}

This can be interpreted as a large deviation estimate
in this context.

\section{Estimates on time} \label {ret random}

Our aim in this section is to estimate the distribution of return times to
$I_n$:
they are concentrated around $c^{-1}_{n-1}$ up to an exponent close to $1$.

The basic estimate is a large deviation estimate which is
proven in the next subsection
(Corollary \ref {large times estimate}) and states that for $k \geq 1$
the set of branches with time larger then
$k c_n^{-4}$ has capacity less then $e^{-k}$.

\subsection{A Large Deviation lemma for times}

Let $r_n(j)$ be such that $R_n|_{I^j_n}=f^{r_n(j)}$.  We will also use the
notation
$r_n(x)=r_n(j^{(n)}(x))$, the $n$-th return time of $x$ (there should be no
confusion for the reader, since we consistently use $j$ for an integer index
and $x$ for a point in the phase space).

Let
$$
A_n(k)=p_{\g_n}(r_n(x) \geq k|I_n)
$$
Since $f$ restricted to the complement of $I_{n+1}$ is hyperbolic, from
Lemma \ref {hyperbol}, it is clear that
$A_n(k)$ decays exponentially with $k$:

\begin{lemma} \label {expan}

With total probability, for all $n>0$, there exists $C>0$, $\lambda>1$
(which depend on $n$) such that $A_n(k)<C \lambda^k$.

\end{lemma}

\begin{pf}

Consider a Markov partition for $f|_{I \setminus I_{n+1}}$, that is, a
finite union of intervals $M_1,...,M_m$ such that

\noindent (1)\, $\cup_{i=1}^m M_i=I \setminus I_{n+1}$,

\noindent (2)\, For every $1 \leq i \leq m$,
$f|_{M_i}$ is a diffeomorphism,

\noindent (3)\,
$f(\cup_{i=1}^m \partial M_i) \subset \cup_{i=1}^m \partial M_i$.

It is easy to see that such a Markov partition also satisfies

\noindent (4)\, For every $1 \leq i \leq m$, either
$$
f(M_i)=\bigcup_{M_j \subset f(M_i)} M_j \quad \text {or} \quad
f(M_i)=I_{n+1} \cup \bigcup_{M_j \subset f(M_i)} M_j.
$$

(To construct such Markov partition, notice first that the boundary of
$I_{n+1}$ is preperiodic to a periodic orbit $q$ (of period $p$).
In particular we have $f^s(\partial I_{n+1})=q$ for some
integer $s>p$.  Let $K$ be
the (finite) set of all $x$ which never enter $\inter I_{n+1}$ and
such that $f^j(x)=q$ for some $j \leq s$.  Since $I_{n+1}$ is nice,
$\partial I_{n+1} \subset K$, and since $s>p$, $q \in K$.  In particular $K$
is forward invariant.  It is easy to see that the connected components of $I
\setminus (K \cup I_{n+1})$ form a Markov partition of $I \setminus
I_{n+1}$.)

It follows that if $f^j(x) \in \cup_{i=1}^m \inter
M_i$, $0 \leq j \leq k$ then there exists a unique interval
$x \in M^k(x)$ such that $f^k|_{M^k(x)}$ is a diffeomorphism onto some
$M_j$.  Notice that if $k \geq 1$, $f(M^k(x))=M^{k-1}(f(x))$.

By Lemma \ref {hyperbol}, if $y \in M^k(x)$, $|Df^k(y)|$ is exponentially
big in $k$.  In particular, $\sum_{j=0}^{k-1} |f^j(M^k(x))|<C'$ for some
constant $C'>0$ independent of $M^k(x)$.  Since $f$ is $C^2$,
$\dist(f|_{M^k(x)})$ is uniformly bounded in $k$.
Notice that the bounds on distortion depend on $n$.
(An alternative to this classical argument is to obtain the
bounded distortion from the negative Schwarzian derivative).

By Lemma \ref {hyperbol} again,
the set of points $x \in I$ which never enter $I_{n+1}$ has
empty interior: for every $T \subset I$ there is an iterate
$f^r(T)$ which intersects $I_{n+1}$ (otherwise the exponentially growing
intervals $f^r(T) \subset I$ would eventually become bigger than $I$).
So there exists $r>0$ such that,
for every $M_j$, there exists $x \in M_j$ and $t_j<r$ with
$f^{t_j}(x) \in \inter I_{n+1}$.
It follows that there exists an interval
$E_j \subset M_j$ such that $f^{t_j}(E_j) \subset \inter I_{n+1}$.

Fixing some $M^k(x)$ with $f^k(M^k(x))=M_j$,
let $E^k(x)=(f^k|_{M^k(x)})^{-1}(E_j)$.  By bounded
distortion, it follows that $|E^k(x)|/|M^k(x)|$ is uniformly bounded from
below independently of $M^k(x)$.  In particular, $p_{2\g} (M^k(x)
\setminus E^k(x)|M^k(x))<\lambda$ for some constant $\lambda<1$.

Let $M^k$ be the union of the $M^k(x)$ and $E^k$ be the union of the
$E^k(x)$.  Then $M^{k+r} \cap E^k=\emptyset$.
In particular, $p_{2\g} (M^{(k+1)r}|I) <
\lambda p_{2\g} (M^{kr}|I)$.

We conclude that
$p_{2\g} (M^k|I_n)<C \lambda^{k/r}$ for some constant $C>0$.
If $k>r_n(0)$, then $M^k \cap I_n$ contains the set of points $x \in I_n$
such that $f^j(x) \notin I_n$, $1 \leq j \leq k$, that is, all points $x \in
I_n$ with $r_n(x)>k$.  Adjusting $C$ and $\lambda$
if necessary, we have $A_n(k)<C \lambda^k$.
\end{pf}

\begin{rem}

It turns out that $\lambda$ depends strongly on $n$.  Indeed, it is possible
to show that $\lambda$ is torrentially close to $1$.  The argument above
does not give any estimate on the behavior of $\lambda$ as $n$ grows,
but it will be used below as the basis of an inductive argument which will
give explicit estimates on $\lambda$ for $n$ big.

\end{rem}

Let $\zeta_n$ be the maximum $\zeta \leq c_{n-1}$ such that for all
$k \geq \zeta^{-1}$ we have
\be
A_n(k) \leq e^{-\zeta k}
\ee
and finally let $\alpha_n=\min_{1 \leq m \leq n} \zeta_m$.

Our main result in this section is to estimate $\alpha_n$.  We will show
that with total probability, for $n$ big we have $\alpha_{n+1} \geq c^4_n$.
For this we will have to do a simultaneous estimate for landing times,
which we define now.

Let $l_n(\d)$ be such that $L_n|_{I^\d_n}=f^{l_n(\d)}$.  We will also use
the notation $l_n(x)=l_n(\d^{(n)}(x))$.

Let us define
\be
B_n(k)=p_{\tg_n}(l_n(x)>k|I_n).
\ee
\be
B^{\tau_n}_n(k)=p_{\tg_n}(l_n(x)>k+r_n(\tau_n)|I^{\tau_n}_n).
\ee


\begin{lemma} \label {landing times}

If $k>c_n^{-3/2} \alpha_n^{-3/2}$ then
\begin{equation} \label {ltimes 1}
B_n(k)<e^{-c_n^{3/2} \alpha_n^{3/2} k},
\end{equation}
and
\begin{equation} \label {ltimes 2}
B^{\tau_n}_n<e^{-c_n^{-3/2} \alpha_n^{3/2}k}.
\end{equation}

\end{lemma}

\begin{pf}

Let us first show (\ref {ltimes 1}).
Let $k>c_n^{-3/2} \alpha_n^{-3/2}$ be fixed.  Let $m_0=\alpha_n^{3/2} k$.

Notice that by Lemma \ref {estimate on m}
\be
p_{\tg_n}(|\d^{(n)}(x)| \geq m_0|I_n) \leq
e^{-c_n^{5/4} \alpha_n^{3/2} k}.
\ee

Fix now $m<m_0$.  Let us estimate
\be
p_{\tg_n} (|\d^{(n)}(x)|=m,l_n(x)>k|I_n).
\ee

For each $\d=(j_1,...,j_m)$ we can associate a sequence of $m$
positive integers $r_i$ such that $r_i \leq r_n(j_i)$ and $\sum r_i=k$.
The average value of $r_i$ is at least $k/m$ so we conclude that
\be
\sum_{r_i \geq k/2m} r_i>k/2.
\ee
Recall also that
\be
\frac {k} {2m}>\frac {1} {(2\alpha_n^{3/2})}>\alpha_n^{-1}.
\ee

Given a sequence of $m$ positive integers $r_i$ as above we
can do the following estimate using Lemma \ref {distortion product}
\begin{align} \label {ri}
p_{\tg_n} (\d^{(n)}(x)=(j_1,...,j_m)&,r_n(j_i) \geq r_i|I_n)\\
\nonumber
&\leq 2^{mn} \prod_{j=1}^m p_{\g_n}(r_n(x) \geq r_j|I_n)\\
\nonumber
&\leq
2^{mn} \prod_{r_j \geq \alpha_n^{-1}} p_{\g_n}(r_n(x) \geq r_j|I_n)\\
\nonumber
&\leq 2^{mn} \prod_{r_j \geq k/2m} e^{-\alpha_n r_j}\\
\nonumber
&\leq 2^{mn} e^{-\alpha_n k/2}.
\end{align}

The number of sequences of $m$ positive integers $r_i$ with sum $k$ is
\begin{align}
\binom {k+m-1} {m-1} &\leq \frac {1}{(m-1)!}
(k+m-1)^{m-1}\\
\nonumber
&\leq \frac{1}{m!}(k+m)^m
\leq \left (\frac {2ek} {m} \right )^m.
\end{align}

Notice that
\begin{align}
2^{mn}\left (\frac {2ek} {m} \right)^m &\leq
\left (\frac {2^{n+3}k} {m} \right )^{\frac {m} {k2^{n+3}} k2^{n+3}}\\
\nonumber
&\leq \left (\frac {2^{n+3}k} {m_0} \right)^{\frac {m_0} {k2^{n+3}} k2^{n+3}}
&&\text {(since $x^{1/x}$ is decreases for $x>e$)}\\
\nonumber
&\leq \left (\frac {2^{n+3}} {\alpha_n^{3/2}} \right )^{m_0} \leq
e^{\alpha_n^{5/4}k}.
\end{align}

So we can finally estimate
\begin{align}
p_{\tg_n} (|\d^{(n)}(x)|=m,l_n(x) \geq k|I_n) &\leq 2^{mn}
\left ( \frac {2ek} {m} \right )^m
e^{-\alpha_n k/2}\\
\nonumber
&< e^{(\alpha_n^{1/4}-1/2)\alpha_n k}.
\end{align}

Summing up on $m$ we get
\begin{align}
p_{\tg_n} (|\d^{(n)}(x)|<m_0&,l_n(x) \geq k|I_n)\\
\nonumber
&\leq m_0 e^{(\alpha_n^{1/4}-1/2)\alpha_n k}\\
\nonumber
&<e^{(2 \alpha_n^{1/4}-1/2) \alpha_n k}
&&(\text {since}
\frac {\ln(m_0)} {k} \leq \frac {\ln(k)} {k} \leq
\alpha_n^{5/4})\\
\nonumber
&\leq e^{-\alpha_n k/3}.
\end{align}

As a direct consequence we get 
\be
B_n(k)<e^{-\alpha_n k/3}+
e^{-c_n^{5/4} \alpha_n^{3/2} k}<e^{-c_n^{3/2} \alpha_n^{3/2} k},
\ee
concluding the proof of (\ref {ltimes 1}).

For the proof of (\ref {ltimes 2}) one proceeds analogously.  Take $k$ and
$m_0$ as before.  By Lemma \ref {8.2} one gets
\be \label {lalala1}
p_{\tg_n}(|\d^{(n)}(x)| \geq m_0|I^{\tau_n}_n) \leq
e^{-c_n^{5/4} \alpha_n^{3/2} k}.
\ee
For any $m<m_0$, if $\d=(\tau_n,j_1,...,j_m)$ and $l_n(\d)>k+r_n(\tau_n)$
then there exists $r_i \leq r_n(j_i)$ with $\sum_{i=1}^m r_i=k$.
Repeating the argument of (\ref {ri})
(and using Lemma \ref {distortion product2} instead of Lemma \ref
{distortion product}) one gets, for any such sequence $r_1,...,r_m$,
\be
p_{\tg_n} (\d^{(n)}(x)=(\tau_n,j_1,...,j_m),r_n(j_i) \geq r_i|I^{\tau_n}_n)
\leq 2^{mn} e^{-\alpha_n k/2}.
\ee
The previous combinatorial estimate can be applied again to obtain
\be \label {lalala2}
p_{\tg_n} (|\d^{(n)}(x)|=m+1,l_n(x)>k+r_n(\tau_n)|I^{\tau_n}_n) <
e^{(\alpha_n^{1/4}-1/2)\alpha_n k}.
\ee
Summing up (\ref {lalala2}) on $m<m_0$ and using (\ref {lalala1}) we
obtain estimate (\ref {ltimes 2}).
\end{pf}

Let $v_n=r_n(0)$ be the return time of the critical point.


\begin{lemma} \label {estimate on v_n}

With total probability, for $n$ large enough,
$$
v_{n+1}<c^{-2}_n \alpha_n^{-2}/2.
$$

\end{lemma}

\begin{pf}

By the definition of $\alpha_n$ and PhPa2, it follows that with total
probability, for $n$ large enough,
$$
r_n(\tau_n)<c_{n-1}^{-1} \alpha_n^{-1}.
$$

Recall that $\d^{(n)}(0)$ is such that $R_n(0) \in C^{\d^{(n)}(0)}_n$.
Using Lemma \ref {landing times}, more precisely estimate (\ref {ltimes 2}),
together with PhPa1, we get with total probability, for $n$ large enough,
$$
l_n(\d^{(n)}(0))-r_n(\tau_n)<n \alpha_n^{-3/2} c_n^{-3/2},
$$
and thus
$$
v_{n+1}<v_n+c_{n-1}^{-1} \alpha_n^{-1}+n \alpha_n^{-3/2}
c_n^{-3/2}<v_n+\alpha_n^{-2}c_n^{-2}/4.
$$
Notice that $\alpha_n$ decreases monotonically, thus
for $n_0$ big enough and for $n>n_0$,
$$
v_{n+1}<v_{n_0}+\sum_{k=n_0}^n \alpha_k^{-2}
c_k^{-2}/4<v_{n_0}+\alpha_n^{-2} c_n^{-2}/3.
$$
which for $n$ big enough implies $v_{n+1}<c_n^{-2} \alpha_n^{-2}/2$.
\end{pf}


\begin{lemma}

With total probability, for $n$ large enough,
$$
\alpha_{n+1} \geq \min \{\alpha_n^4, c_n^4\}.
$$

\end{lemma}

\begin{pf}

Let $k \geq \max\{\alpha_n^{-4},c_n^{-4}\}$.
From Lemma \ref {estimate on v_n} one
immediately sees that if $r_{n+1}(j) \geq k$ then $R_n(I^j_{n+1}))$ is
contained on some $C^\d_n$ with $l_n(\d)=r_{n+1}(j)-v_n \geq
k/2 \geq n \alpha_n^{-3/2} c_n^{-3/2}$.

Applying Lemma \ref {landing times} we have
$B_n(k/2)<e^{-\alpha_n^{3/2} c_n^{3/2} k/2}$.

Applying Remark \ref {remark 2} we get
$$
A_{n+1}(k)<e^{-k \alpha_n^{3/2}
c_n^{3/2}/200}<e^{-k \min\{\alpha_n^{4},c_n^{4}\}}.
$$
\end{pf}

Since $c_n$ decreases torrentially, we get


\begin{cor} \label {large times estimate}

With total probability, for $n$ large enough $\alpha_{n+1} \geq c_n^4$.

\end{cor}

\begin{rem} \label {v_n}

In particular, using Lemma \ref {estimate on v_n}, for $n$ big,
$v_n<c_{n-1}^{-2} \alpha_{n-1}^{-2}/2<c_{n-1}^{-4}$.

\end{rem}

\subsection{Consequences}

\begin{lemma} \label {l_n}

With total probability, for all $n$ sufficiently large we have
\begin{align}
\label{d1}
&p_{\tg_n}(l_n(x)<c_n^{-1+\epsilon}|I_n)<
c^{\epsilon/2}_n,\\
\label{d2}
&p_{\tg_n} (l_n(x)>c_n^{-1-5\epsilon/3}|I_n) \leq
e^{-c^{-\epsilon/4}_n},\\
\label{d3}
&p_{\tg_n}(l_n(x)-r_n(x)<c_n^{-1+\epsilon}|I^{\tau_n}_n)<
c^{\epsilon/2}_n,\\
\label{d4}
&p_{\tg_n} (l_n(x)-r_n(x)>c_n^{-1-5\epsilon/3}|I^{\tau_n}_n)
\leq e^{-c^{-\epsilon/4}_n}.
\end{align}

\end{lemma}

\begin{pf}

We will concentrate on estimates (\ref {d1}) and (\ref {d2}), since (\ref
{d3}) and (\ref {d4}) are analogous.

We have $l_n(\d) \geq |\d|$, and
from Lemma \ref {estimate on m}
$$
p_{\tg_n}(|\d^{(n)}(x)| \leq
c^{-1+\epsilon}_n|I_n) \leq
c_n^{\epsilon/2},
$$
which implies (\ref {d1}).

On the other hand, by the same lemma,
$$
p_{\tg_n}(|\d^{(n)}(x)|
\geq c^{-1-\epsilon}_n|I_n) \leq
e^{-c^{-\epsilon/2}_n}.
$$
Defining
$$
X_m=\bigcup_{
\ntop
{\d=(j_1,...,j_m),}
{r_n(j_m)>c_n^{-\epsilon/2} c_{n-1}^{-4}}
}
I^{\d}_n
$$
we have
$$
p_{\tg_n}(X_m|I_n) \leq
2^n e^{-c^{-\epsilon/2}_n}<e^{-c^{-\epsilon/3}_n}.
$$

Since 
$$
c^{-1-\epsilon}_n c^{-\epsilon/2}_n c^{-4}_{n-1}<c^{-1-5\epsilon/3}_n,
$$
we conclude that if $x$ satisfies $l_n(x)>c_n^{-1-5\epsilon/3}$ and
$|\d_n(x)|<c_n^{-1-\epsilon}$ then $x$ belongs to some
$X_m$ with $1 \leq m \leq c_n^{-1-\epsilon}$.
So we get
$$
p_{\tg_n} (l_n(x)>c_n^{-1-5\epsilon/3}|I_n) \leq
e^{-c^{-\epsilon/2}_n}+c_n^{-1-\epsilon}
e^{-c_n^{-\epsilon/3}}<e^{-c_n^{-\epsilon/4}}
$$
which implies (\ref {d2}).
\end{pf}


\begin{cor} \label {precise time estimate}

With total probability, for all $n$
sufficiently large we have
\begin{align}
\label{d5}
&p_{\g_{n+1}}
(r_{n+1}(x)<c_n^{-1+\epsilon}|I_{n+1})<c^{\epsilon/10}_n,\\
\label{d6}
&p_{\g_{n+1}} (r_{n+1}(x)>c_n^{-1-2\epsilon}|I_{n+1})
\leq e^{-c^{-\epsilon/5}_n} \leq c^n_n.
\end{align}

\end{cor}

\begin{pf}

Notice that $r_{n+1}(j)=v_n+l_n(\d)$, where $R_n(I^j_{n+1}) \subset C^\d_n$.
By Remark \ref {v_n}, we can estimate $v_n<c_{n-1}^{-10}$.
The distribution of $r_{n+1}(j)-v_n$ can be then estimated
by the distribution of $l_n(\d)$ from Lemma \ref {l_n},
with a slight loss given by
Remark \ref {remark 2}.
\end{pf}

Using PhPa2 we get


\begin{lemma} \label {first iterate}

With total probability, for all $n$ sufficiently big
\begin{equation}
\lim_{n \to \infty} \frac {\ln(r_n(\tau_n))} {\ln(c^{-1}_{n-1})}=1.
\end{equation}

\end{lemma}

\begin{cor} \label {l_n 2}

With total probability, for all $n$ sufficiently large we have
\begin{align}
&p_{\tg_n}(l_n(x)<c_n^{-1+\epsilon}|I^{\tau_n}_n) \leq
c^{\epsilon/10}_n,\\
&p_{\tg_n} (l_n(x)>c_n^{-1-11\epsilon/6}|I^{\tau_n}_n) \leq
e^{-c^{-\epsilon/5}_n}.
\end{align}

\end{cor}

\begin{pf}

Just use Lemma \ref {first iterate} together with estimates (\ref {d3}) and
(\ref {d4}) of Lemma \ref {l_n}.
\end{pf}


\begin{cor} \label {growth of v_n}

With total probability,
$$
\lim_{n \to \infty} \frac {\ln(v_{n+1})} {\ln(c^{-1}_n)}=1.
$$

\end{cor}

\begin{pf}

Notice that $v_{n+1}=v_n+l_n(\d)$ where
$R_n(0) \in C^\d_n$.
Using Corollary \ref {l_n 2} and PhPa1 we get
$c_n^{-1+\epsilon}<l_n(\d)<c_n^{-1-11\epsilon/6}$.  By Remark \ref {v_n},
$v_n<c_{n-1}^{-10}$, so
$c_n^{-1+\epsilon}<v_{n+1}<c_n^{-1-2\epsilon}$.  Letting $\epsilon$ go to
$0$ we get the result.
\end{pf}


\begin{rem} \label {precise rate of c_n}

Using Lemma \ref {away from the boundary}, we see that
$|R_n(I_{n+1})|>2^{-n}|I_n|$.  Since $|Df(x)|<4$, $x \in I$,
it follows that $|DR_n(x)|<4^{v_n}$, $x \in I_{n+1}$.
In particular, Corollary \ref {growth of v_n} implies that
with total probability, for all $\epsilon>0$,
for all $n$ big enough,
$$
2^{-n} c_n^{-1}<\frac {|R_n(I_{n+1})|}
{|I_{n+1}|}<4^{v_n}<4^{c_{n-1}^{-1-\epsilon}},
$$
so $\ln (c_n^{-1})<c_{n-1}^{-1-2\epsilon}$.
This implies together with Corollary \ref {c_n
torrential} that
$$
\lim_{n \to \infty} \frac {\ln(\ln(c_n^{-1}))} {\ln(c_{n-1}^{-1})}=1,
$$
so $c^{-1}_n$ grows torrentially (and not faster).

\end{rem}

\section{Dealing with hyperbolicity}

In this section we show by an inductive process that the great
majority of branches are reasonably hyperbolic.
In order to do that, in the
following subsection, we define some classes of branches with
``very good'' distribution of times and which
are not too close to the critical point.  The definition of
very good distribution of times has an inductive
component: they are composition of many very good branches of the previous
level.  The fact that most branches are very good is related to the
validity of some kind of Law of Large Numbers estimate.

\subsection{Some kinds of branches and landings}

\subsubsection{Standard landings}

Let us define the set of standard landings of level $n$,
$LS(n) \subset \Omega$ as the set
of all $\d=(j_1,...,j_m)$ satisfying the following.

\begin{description}

\item[LS1] ($m$ not too small or large) $c^{-1/2}_n<m<c^{-1-2\epsilon}_n$,

\item[LS2] (No very large times) $r_n(j_i)<c^{-14}_{n-1}$ for all $i$.

\item[LS3] (Short times are sparse in large enough initial segments)
For $c^{-2}_{n-1} \leq k \leq m$
$$
\#\{1 \leq i \leq k,\, r_n(j_i)<c^{-1+2\epsilon}_{n-1}\} <
(6 \cdot 2^n) c^{\epsilon/10}_{n-1} k,
$$
  
\item[LS4] (Large times are sparse in large enough initial segments)
For $c^{-1/n}_n \leq k \leq m$
$$
\#\{1 \leq i \leq k,\, r_n(j_i)>c^{-1-2\epsilon}_{n-1}\} <
(6 \cdot 2^n) e^{-c^{\epsilon/5}_{n-1}} k.
$$

\end{description}


\begin{lemma} \label {standard landing}

With total probability, for all $n$ sufficiently big we have
\begin{equation} \label {k1}
p_{\tg_n}(\d^{(n)}(x) \notin LS(n)|I_n)<c_n^{1/3}/2,
\end{equation}
\begin{equation} \label {k2}
p_{\tg_n}(\d^{(n)}(x) \notin LS(n)|I^{\tau_n}_n)<c_n^{1/3}/2.
\end{equation}

\end{lemma}

\begin{pf}

Let us start with estimate (\ref {k1}) (on $I_n$).  Let us estimate
the complement of the set of
landings which violate each item of the definition.

(LS1)\, This was already estimated before (see
Lemma \ref {estimate on m}),
an upper bound is $c_n^{1/3}/3$ (using $\epsilon$ small).

(LS2)\, By Corollary \ref {large times estimate}
the $\g_n$-capacity of $\{r_n(x)>c_{n-1}^{-14}\}$ is at most
$e^{-c^{-10}_{n-1}} \ll c_n^3$.  Using Lemma \ref {2n},
we see that the $\tg_n$-capacity of the set of
$\d=(j_1,...,j_m)$ with $r_n(j_i)>c_n^{-14}$
for some $i \leq c_n^{-1-2\epsilon}$ (in particular for some $i \leq m$
if $m$ is as in LS1)
is bounded by $2^n c_n^{-1-2\epsilon} c_n^3 \ll c_n$.

(LS3)\, This is a large deviation estimate, so we follow the ideas of \S
\ref {more tree}, particularly estimate (\ref {q}).
Put $q=c_{n-1}^{\epsilon/10}$.  By estimate (\ref
{d5}) of Corollary \ref {precise time estimate},
we can estimate the $\tilde \gamma_n$-capacity corresponding to the
violation of LS3 for some fixed $c_{n-1}^{-2}
\leq k \leq c_n^{-1-2\epsilon}$ by
$$
\left (\frac {1} {2} \right )^{(6 \cdot 2^n) q k} \leq
\left (\frac {1} {2} \right)^{c^{-3/2}_{n-1}} \ll c_n^3.
$$
Summing up over $k \leq c_n^{-1-2\epsilon}$
(and in particular for $k \leq m$ as in LS1) we get the upper
bound $c_n$.

(LS4)\, We use the same method of the previous item.  Put
$q=e^{-c_{n-1}^{-\epsilon/5}}$.  By estimate (\ref {d6}) of
Corollary \ref {precise time estimate}, we can bound
the $\tilde \gamma_n$-capacity corresponding to the violation of LS4
for some fixed $c_n^{-1/n} \leq k \leq c_n^{-1-2\epsilon}$ by
$$
\left (\frac {1} {2} \right )^{(6 \cdot 2^n) q k} \ll c_n^3.
$$
Summing up over $k \leq c_n^{-1-2\epsilon}$
(and in particular for $k \leq m$ as in LS1) we get the upper
bound $c_n$.

Adding the losses of the four items, we get the estimate (\ref {k1}).
To get estimate (\ref {k2}) (on $I^{\tau_n}_n$),
the only item which changes is LS2: we have to be careful since if
$r_n(\tau_n)$ is very large then automatically LS2 is violated for every
$\d$ which starts by $\tau_n$.  But this was taken care in Lemma \ref {first
iterate}, and with this observation the estimates are the same.
\end{pf}

\subsubsection{Very good returns and excellent landings}

Define the set of very good returns,
$VG(n_0,n) \subset \Z \setminus \{0\}$, $n_0,n \in \N$,
$n \geq n_0$ by induction as follows.  We let
$VG(n_0,n_0)=\Z \setminus \{0\}$ and supposing
$VG(n_0,n)$ already defined, let us define $LE(n_0,n) \subset LS(n)$
(excellent landings) as the set of standard landings satisfying the
following extra condition:

\begin{description}

\item[LE] (Not very good moments are sparse in large enough initial segments)
For all $c^{-2}_{n-1}<k \leq m$
$$
\#\{1 \leq i \leq k,\, j_i \notin VG(n_0,n)\} <
(6 \cdot 2^n) c^{1/20}_{n-1} k,
$$

\end{description}

And we define $VG(n_0,n+1)$ as the set of $j$ such that
$R_n(I^j_{n+1})=C^{\d}_n$ with $\d \in LE(n_0,n)$ and 
the satisfying the extra condition:

\begin{description}

\item[VG] (distant from $0$) The distance of $I^j_{n+1}$ to $0$ is bigger
than $c_n^{1/3}|I_{n+1}|$.

\end{description}


\begin{lemma} \label {induction step}

With total probability, for all $n_0$ sufficiently big and all
$n \geq n_0$, if
\begin{equation}
p_{\g_n}(j^{(n)}(x) \notin VG(n_0,n)|I_n)<c_{n-1}^{1/20}
\end{equation}
then
\begin{align}
\label{le1}
&p_{\tg_n}(\d^{(n)}(x) \notin LE(n_0,n)|I_n)<c_n^{1/3},\\
\label{le2}
&p_{\tg_n}(\d^{(n)}(x) \notin LE(n_0,n)|I^{\tau_n}_n)<c_n^{1/3}.
\end{align}

\end{lemma}

\begin{pf}

We first use Lemma \ref {standard landing}
to estimate the $\tg_n$-capacity of branches not in $LS(n)$ by
$c_n^{1/3}/2$.

Let $q=c^{1/20}_{n-1}$.
Using the hypothesis and estimate (\ref {q}) of \S \ref {more tree}
(see also the estimate of the complement of
LS3 in Lemma \ref {standard landing})
we first estimate
the $\tg_n$-capacity of the set of landings which violate LE
for a specific value of $k$ with $k \geq c_{n-1}^{-2}$
by $(1/2)^{(6 \cdot 2^n) qk}$
and then summing up on $k$ we get
$$
\sum_{k \geq
c^{-2}_{n-1}} \left (\frac {1} {2} \right)^{(6 \cdot 2^n) q k} \ll c_n.
$$

This argument works both for (\ref {le1}) (in $I_n$) and (\ref {le2})
(in $I^{\tau_n}_n$).
\end{pf}


\begin{lemma} \label {very good capacity}

With total probability, for all $n_0$ sufficiently big and for all $n \geq
n_0$,
\begin{equation} \label {vg}
p_{\g_n}(j^{(n)}(x) \notin VG(n_0,n)|I_n)<c_{n-1}^{1/20}.
\end{equation}

\end{lemma}

\begin{pf}

It is clear that with total probability, for $n_0$ sufficiently big and
$n \geq n_0$, the set of branches $I^j_n$ at distance at least
$c_{n-1}^{1/3}|I_n|$ of $0$
has $\g_n$-capacity bounded by $c_{n-1}^{1/8}$.

For $n=n_0$, (\ref {vg}) holds (since all branches
are very good except the central one).  Using
Lemma \ref {induction step}, if (\ref {vg}) holds for $n$ then (\ref {le1})
also holds for $n$.  Pulling back estimate (\ref {le1}) by $R_n|_{I_{n+1}}$
(using Remark \ref {remark 2}), we get (\ref {vg}) for $n+1$.  By induction
on $n$, (\ref {vg}) holds for all $n \geq n_0$.
\end{pf}

Using PhPa2 we get (using the measure-theoretical argument of Lemma \ref
{simp})


\begin{lemma} \label {lands very good}

With total probability, for all $n_0$ big enough, for all $n$ big enough,
$\tau_n \in VG(n_0,n)$.

\end{lemma}


\begin{lemma} \label {very good return time}

With total probability, for all $n_0$ big enough and for all $n \geq n_0$,
if $j \in VG(n_0,n+1)$ then
$$
\frac{1}{2} m c^{-1+2 \epsilon}_{n-1}<r_{n+1}(j)<
2 m c^{-1-2 \epsilon}_{n-1},
$$
where as usual, $m$ is such that $R_n(I^j_{n+1})=C^{\d}_n$ and
$\d=(j_1,...,j_m)$.

\end{lemma}

\begin{pf}

Notice that $r_{n+1}(j)=v_n+\sum r_n(j_i)$.
To estimate the total time $r_{n+1}(j)$
from below we use LS3 and get
$$
\frac{1}{2} m c^{-1+2 \epsilon}_{n-1}<(1-6 \cdot 2^n c^{\epsilon/10}_n) m
c^{-1+2 \epsilon}_{n-1} < r_{n+1}(j).
$$

To estimate from above, we notice that $v_n<c^{-4}_{n-1}$ and by LS2 and LS4
$$
\sum_{r_n(j_i)>c^{-1-2 \epsilon}_{n-1}} r_n(j_i)<
6 \cdot 2^n c^{-14}_{n-1} e^{-c^{-\epsilon/5}_{n-1}} m < m,
$$
so
$$
r_{n+1}(j)<m c^{-1-2 \epsilon}_{n-1}+m+c^{-4}_{n-1}<2 m
c^{-1-2 \epsilon}_{n-1}.
$$
\end{pf}

\begin{rem} \label {very good time again}

Using LS1 we get the estimate $c_n^{-1/2}<r_{n+1}(j)<c_n^{-1-3\epsilon}$ for
$j \in VG(n_0,n+1)$.

\end{rem}

Let $j \in VG(n_0,n+1)$.  We can write
$R_{n+1}|_{I^j_{n+1}}=f^{r_{n+1}(j)}$, that is, a big iterate of $f$.  One
may consider which proportion of those iterates belong to very good
branches of the previous level.  More generally, we can truncate the return
$R_{n+1}$, that is, we may consider
$k<r_{n+1}(j)$ and ask which proportion of iterates up to $k$ belong to very
good branches.


\begin{lemma} \label {very good partial time}

With total probability, for all $n_0$ big enough and for all $n \geq n_0$,
the following holds.

Let $j \in VG(n_0,n+1)$, and let $\d=(j_1,...,j_m)$ be such that
$R_n(I^j_{n+1})=C^\d_n$.  Let $m_k$ be biggest possible with
$$
v_n+\sum_{j=1}^{m_k} r_n(j_i) \leq k
$$
(the amount of full returns to level $n$ before time $k$) and let
$$
\beta_k=\sum_{\ntop {1 \leq i \leq m_k,}
{j_i \in VG(n_0,n)}} r_n(j_i)
$$
(the total time spent in full returns to level $n$
which are very good before time $k$).
Then $1-\beta_k/k<c_{n-1}^{1/100}$ if $k>c_n^{-2/n}$.

\end{lemma}

\begin{pf}

Let us estimate first the time $i_k$ which is not spent on non-critical
full returns:
$$
i_k=k-\sum_{j=1}^{m_k} r_n(j_i).
$$
This corresponds exactly to $v_n$ plus some incomplete part of the return
$j_{m_{k+1}}$.  This part can be bounded by $c_{n-1}^{-4}+c_{n-1}^{-14}$
(use Corollary \ref {growth of v_n} to estimate $v_n$ and LS2 to estimate
the incomplete part).

Using LS2 we conclude now that
$$
m_k>(k-c_{n-1}^{-4}-c_{n-1}^{-14})c_{n-1}^{14}>c_n^{-1/n}
$$
so $m_k$ is not too small.

Let us now estimate the contribution $h_k$ from full returns $j_i$
with time higher than $c_{n-1}^{-1-2\epsilon}$.
Since $m_k$ is big, we can use
LS4 to conclude that the number of such high time returns
must be less than $c_{n-1}^n m_k$, so their total time is at most
$c_{n-1}^{n-14} m_k$.

The non very good full returns on the other hand can be estimated by LE
(given the estimate on $m_k$), they are at most $c_{n-1}^{1/21} m_k$. 
So we can estimate the total time $l_k$ of non very good full returns
with time less then $c_{n-1}^{-1-2\epsilon}$ by
$$
c_{n-1}^{1/25}c_{n-1}^{-1-2\epsilon} m_k.
$$
Since $m_k$ is big, we can use LS3 to estimate the proportion of branches
with not too small time, so we conclude that at most
$c_{n-1}^{\epsilon/11} m_k$
branches are not very good or have time less than
$c_{n-1}^{-1+2\epsilon}$, so $\beta_k$ can be estimated from below as
$$
(1-c^{\epsilon/11}_{n-1}) c^{-1+2\epsilon}_{n-1} m_k.
$$

It is easy to see then that $i_k/\beta_k \ll c^{1/100}_{n-1}$,
$h_k/\beta_k \ll c^{1/100}_{n-1}$.
If $\epsilon$ is small enough, we also have
$$
l_k/\beta_k<2 c_{n-1}^{1/25-4 \epsilon} < c^{1/90}_{n-1}.
$$
So $(i_k+h_k+l_k)/\beta_k$ is less then $c^{1/100}_{n-1}$.
Since $i_k+h_k+l_k+\beta_k=k$ we have $1-\beta_k/k<(i_k+h_k+l_k)/\beta_k$.
\end{pf}

\subsubsection{Cool landings}

Let us define the set of cool landings
$LC(n_0,n) \subset \Omega$, $n_0,n \in \N$,
$n \geq n_0$ as the set of all $\d=(j_1,...,j_m)$ in $LE(n_0,n)$
satisfying

\begin{description}

\item[LC1] (Starts very good) $j_i \in VG(n_0,n)$, $1 \leq i \leq
c_{n-1}^{-1/30}$,

\item[LC2] (Short times are sparse in large enough initial segments)
For $c^{-\epsilon/5}_{n-1} \leq k \leq m$
$$
\#\{1 \leq i \leq k,\, r_n(j_i)<c^{-1+2\epsilon}_{n-1}\} <
(6 \cdot 2^n) c^{\epsilon/10}_{n-1} k,
$$

\item[LC3] (Not very good moments are sparse in large enough initial segments)
For all $c^{-1/30}_{n-1} \leq k \leq m$
$$
\#\{1 \leq i \leq k,\, j_i \notin VG(n_0,n)\} <
(6 \cdot 2^n) c^{1/60}_{n-1} k,
$$

\item[LC4] (Large times are sparse in large enough initial segments)
For $c_{n-1}^{-200} \leq k \leq m$
$$
\#\{1 \leq i \leq k,\, r_n(j_i)>c^{-1-2\epsilon}_{n-1}\} <
(6 \cdot 2^n) c^{100}_{n-1} k,
$$

\item[LC5] (Starts with no large times) $r_n(j_i)<c_{n-1}^{-1-2\epsilon}$,
$1 \leq i \leq e^{c_{n-1}^{-\epsilon/5}/2}$.

\end{description}

Notice that LC4 and LC5 overlap, since
$c_{n-1}^{-200}<e^{c_{n-1}^{-\epsilon/5}/2}$ as do LC1 and LC3.
From this we can conclude that we can control the proportion of
large times or non very
good times in all moments (and not only for large enough initial segments).


\begin{lemma}

With total probability, for all $n_0$ sufficiently big and all $n \geq n_0$,
\begin{equation} \label {k11}
p_{\tg_n}(\d^{(n)}(x) \notin LC(n_0,n)|I_n)<c_{n-1}^{1/100}
\end{equation}
and for all $n$ big enough
\begin{equation} \label {k12}
p_{\tg_n}(\d^{(n)}(x) \notin LC(n_0,n)|I^{\tau_n}_n)<c_{n-1}^{1/100}.
\end{equation}

\end{lemma}

\begin{pf}

\comm{
The estimates follow the ideas of the proof of Lemma \ref {standard
landing}.  More precisely, LC1 and LC5 can be estimated as we did for
LS2, while the remaining are large deviation estimates where we can
argue as for LS3.  The computations below indicate what is lost
going from excellent to cool:
}

We follow the ideas of the proof of Lemma \ref {standard
landing}.  Let us start with estimate (\ref {k11}).
Notice that by Lemmas \ref {very good capacity} and \ref {induction step}
we can estimate the
$\tg_n$-capacity of the complement of excellent landings by $c_n^{1/3}$.
The computations below indicate what is lost
going from excellent to cool due to each item of the definition:

(LC1)\, This is a direct estimate analogous to LS2.
By Lemma \ref {very good capacity},
the $\g_n$-capacity of the complement of very good branches is
bounded by $c_{n-1}^{1/20}$, so an upper bound for the
$\tg_n$-capacity of the set of landings
which do not start with $c_{n-1}^{-1/30}$
very good branches is given by
$$
2^n c_{n-1}^{1/20} c_{n-1}^{-1/30} \ll c_{n-1}^{1/100}.
$$

(LC2)\, This is essentially the same
large deviation estimate of LS3.  We put $q=c_{n-1}^{\epsilon/10}$.
By estimate (\ref {d5}) of Corollary \ref {precise time estimate},
the $\tg_n$-capacity of the set of landings
violating LC2 for a specific value of $k$ is bounded by $(1/2)^{(6 \cdot
2^n) qk}$, and summing up on $k$ (see also estimate (\ref {sumq}))
we get the upper bound
$$
\sum_{k \geq c_{n-1}^{-\epsilon/5}}
\left (\frac {1} {2} \right)^{(6 \cdot 2^n) c_{n-1}^{\epsilon/10} k} \leq
(2^{-n} c^{-\epsilon/10}_{n-1}) \left (\frac {1} {2} \right)^{(6 \cdot 2^n)
c_{n-1}^{-\epsilon/10}}
\ll c_{n-1}^{1/100}.
$$

(LC3)\, This is analogous to the previous item, we set
$q=c_{n-1}^{1/60}$ and using Lemma \ref {very good capacity}
we get an upper bound
$$
\sum_{k \geq c_{n-1}^{-1/30}}\left (\frac {1} {2} \right)^{(6 \cdot 2^n)
c_{n-1}^{1/60} k} \leq
(2^{-n} c^{-1/60}_{n-1}) \left (\frac {1} {2} \right)^{(6 \cdot 2^n)
c_{n-1}^{-1/60}} \ll c_{n-1}^{1/100}.
$$

(LC4)\, As before, we set $q=c_{n-1}^{100}$
and using estimate (\ref {d6}) of Corollary \ref {precise time estimate}
we get
$$
\sum_{k \geq c_{n-1}^{-200}}\left (\frac {1} {2} \right)^{(6 \cdot 2^n)
c_{n-1}^{100} k} \leq
(2^{-n} c^{-100}_{n-1}) \left (\frac {1} {2} \right)^{(6 \cdot 2^n)
c_{n-1}^{-100}} \ll c_{n-1}^{1/100}.
$$

(LC5)\, This is a direct estimate as LC1, using estimate (\ref {d6}) of
Corollary \ref {precise time estimate} we get
$$
2^n e^{-c_{n-1}^{\epsilon/5}} e^{c_{n-1}^{-\epsilon/5}/2} \ll
c_{n-1}^{1/100}.
$$

\comm{
\begin{align}
\tag{LC1}&
2^n c_{n-1}^{1/20} c_{n-1}^{-1/30} \ll c_{n-1}^{1/100},\\
\tag{LC2}&
\sum_{k \geq c_{n-1}^{-\epsilon/5}}
\left (\frac {1} {2} \right)^{(6 \cdot 2^n) c_{n-1}^{\epsilon/10} k} \leq
(2^{-n} c^{-\epsilon/10}_{n-1}) \left (\frac {1} {2} \right)^{(6 \cdot 2^n)
c_{n-1}^{-\epsilon/10}}
\ll c_{n-1}^{1/100},\\
\tag{LC3}&
\sum_{k \geq c_{n-1}^{-1/30}}\left (\frac {1} {2} \right)^{(6 \cdot 2^n)
c_{n-1}^{1/60} k} \leq
(2^{-n} c^{-1/60}_{n-1}) \left (\frac {1} {2} \right)^{2^n c_{n-1}^{-1/60}}
\ll c_{n-1}^{1/100},\\
\tag{LC4}&
\sum_{k \geq c_{n-1}^{-200}}\left (\frac {1} {2} \right)^{(6 \cdot 2^n)
c_{n-1}^{100} k} \leq
(2^{-n} c^{-100}_{n-1}) \left (\frac {1} {2} \right)^{2^n c_{n-1}^{-100}}
\ll c_{n-1}^{1/100},\\
\tag{LC5}&
2^n e^{-c_{n-1}^{\epsilon/5}} e^{c_{n-1}^{-\epsilon/5}/2} \ll
c_{n-1}^{1/100}.
\end{align}
}

Putting those together, we obtain (\ref {k11}).
For (\ref {k12}), we must be careful to have $\tau_n \in VG(n_0,n)$ and
$r_n(\tau_n)<c_{n-1}^{-1-2 \epsilon}$, otherwise we would have
immediate problems due to LC1 and LC5.  But we took care of those properties
in Lemmas \ref {lands very good} and \ref {first iterate}, and with this
observation the estimates are the same as before.
\end{pf}

Transferring the result to the parameter, using PhPa1,
we get (using the measure-theoretical argument of Lemma \ref {simp})

\begin{lemma} \label {lands cool}

With total probability, for all $n_0$ big enough, for all $n$ big enough we
have $R_n(0) \in C^\d_n$ with $\d \in LC(n_0,n)$.

\end{lemma}

\subsection{Hyperbolicity} \label {hyp branches}

\subsubsection{Preliminaries}

For $j \neq 0$, we define
$$
\lambda_n(j)=\inf_{x \in I^j_n} \frac
{\ln |DR_n(x)|} {r_n(j)}.
$$
And $\lambda_n=\inf_{j \neq 0} \lambda_n(j)$.  As a consequence of the
exponential estimate on distortion for returns
(which competes with torrential expansion from the decay of geometry),
together with hyperbolicity of $f$ in
the complement of $I^0_n$ we immediately have the following


\begin{lemma}

With total probability, for all $n$ sufficiently big, $\lambda_n>0$.

\end{lemma}

\begin{pf}

By Lemma~\ref {hyperbol}, there exists a
constant $\tilde \lambda_n>0$ such that
each periodic orbit $p$ of $f$ whose orbit is entirely contained in
the complement of $I_{n+1}$ must satisfy $\ln |Df^m(p)|>\tilde \lambda_n m$,
where $m$ is the period of $p$.  On the other hand, each non-central branch
$R_n|_{I^j_n}$ has a fixed point.  By Lemma~\ref {distortion},
$\sup \dist(R_n|_{I^j_n}) \leq 2^n$ and of course $\lim_{j \to \pm \infty}
r_n(j)=\infty$, so we have
$$
\liminf_{j \to \pm \infty} \lambda_n(j) \geq \tilde \lambda_n.
$$

On the other hand, for any $j \neq 0$, $\lambda_n(j)>0$ by Lemma~\ref
{distortion}, so $\lambda_n>0$.
\end{pf}

\subsubsection{Good branches}

The ``minimum hyperbolicity'' $\liminf \lambda_n$ of the parameters we
will obtain will in fact be positive, as it follows from one of the
properties of Collet-Eckmann parameters (uniform hyperbolicity on periodic
orbits, see \cite {NS}), together with our estimates on distortion.

However our strategy is not to show that the minimum hyperbolicity
is positive, but that the typical value of $\lambda_n(j)$
stays bounded away from $0$
as $n$ grows (and is in fact bigger than $\lambda_{n_0}/2$
for $n>n_0$ big).  Since we also have to estimate the
hyperbolicity of truncated branches it will be convenient
to introduce a new class of branches with good hyperbolic properties.

We define the set of good returns
$G(n_0,n) \subset \Z \setminus \{0\}$, $n_0,n \in \N$,
$n \geq n_0$ as the set of all $j$ such that

\begin{description}

\item[G1] (hyperbolic return)
$$
\lambda_n(j) \geq \lambda_{n_0} \frac {1+2^{n_0-n}} {2},
$$

\item[G2] (hyperbolicity in truncated return)
for $c_{n-1}^{-3/(n-1)} \leq k \leq r_n(j)$ we have
$$
\inf_{x \in I^j_n} \frac {\ln |Df^k(x)|} {k} \geq
\lambda_{n_0} \frac {1+2^{n_0-n+1/2}} {2}-c_{n-1}^{2/(n-1)}.
$$

\end{description}

Notice that since $c_n$ decreases torrentially, for $n$ sufficiently big G2
implies that if $j$ is good
then for $c_{n-1}^{-3/(n-1)} \leq k \leq r_n(j)$ we have
$$
\inf_{x \in I^j_n}\frac {\ln |Df^k(x)|} {k} \geq
\lambda_{n_0} \frac {1+2^{n_0-n}} {2}.
$$


\begin{lemma} \label {very good is good}

With total probability, for all $n_0$ big enough and for all $n>n_0$,
$VG(n_0,n) \subset G(n_0,n)$.

\end{lemma}

\begin{pf}

Let us prove that if G1 is satisfied
for all $j \in VG(n_0,n)$, then $VG(n_0,n+1) \subset G(n_0,n+1)$ (notice
that by definition of $\lambda_{n_0}$ the hypothesis
is satisfied for $n_0$).  Fix $j \in VG(n_0,n+1)$ and define
$$
a_k=\inf_{x \in I^j_{n+1}} \frac {\ln |Df^k(x)|} {k},
$$
and let us consider values of $k$ in the range
$c_n^{-3/n} \leq k \leq r_{n+1}(j)$ (notice that if $k=r_{n+1}(j)$ belongs
to this range by Remark \ref {very good time again}).

We let (as usual) $R_n(I^j_{n+1}) \subset C^\d_n$, $\d=(j_1,...,j_m)$.
Notice that by Corollary \ref {growth of v_n}, $v_n<c_{n-1}^{-4}<k$.
Let us say that $j_i$ was completed before $k$ if
$v_n+r_n(j_1)+...+r_n(j_i) \leq k$.  We let the queue be defined as
$$
q_k=\inf_{x \in C^\d_n} \ln |Df^{k-r} \circ f^r(x)|
$$
where
$r=v_n+r_n(j_1)+...+r_n(j_{m_k})$ with $j_{m_k}$ the last complete return.

Let us show first that
$|DR_n(x)|>1$ if $x \in I^j_{n+1}$.
Indeed, by Lemma \ref {decomposition}, $DR_n|_{I_{n+1}}=\phi \circ f$,
where $\phi$ has small distortion, so by Lemma \ref {away from the
boundary},
$$
|D\phi(x)|>\frac {|R_n(I_{n+1})|} {2|f(I_{n+1})|}>
\frac {2^{-n} |I_n|} {|I_{n+1}|^2},
$$
while by VG,
$|Df(x)|=|2 x| \geq c_n^{1/3} |I_{n+1}|$, so $|DR_n(x)|>c_n^{-1/2}$.

By Lemma \ref {distortion},
any complete return before $k$ produces some expansion,
(that is, the absolute value of the derivative of such return is at least
$1$).  On the other hand, $-q_k$ can
be bounded from above by $-\ln(c_n c_{n-1}^5)$
using Lemma \ref {lower bound}.  We have
$$
-\frac {q_k} {k} \leq \frac {-\ln(c_n c_{n-1}^5)}
{c_n^{-3/n}} \ll c_n^{2/n}.
$$

Now we use Lemma \ref {very good partial time} and get
\begin{align*}
a_k &> \frac {\beta_k} {k} \frac {\lambda_{n_0}(1+2^{n_0-n})} {2}-
\frac {-q_k} {k}\\
&\geq \frac {\lambda_{n_0}(1+2^{n_0-n-1/2})} {2}-\frac {-q_k} {k}
\end{align*}
which gives G2.  If $k=r_{n+1}(j)$ then $q_k=0$, which gives G1.
\end{pf}

\subsubsection{Hyperbolicity in cool landings}


\begin{lemma} \label {cool hyperbolicity}

With total probability, if $n_0$ is sufficiently big, for all $n$
sufficiently big, if $\d \in LC(n_0,n)$ then for all
$c_{n-1}^{-4/(n-1)}<k \leq l_n(\d)$,
$$
\inf_{x \in C^\d_n} \frac {\ln |Df^k(x)|} {k} \geq
\frac {\lambda_{n_0}} {2}.
$$

\end{lemma}

\begin{pf}

Fix such $\d \in LC(n_0,n)$, and let as usual $\d=(j_1,...,j_m)$.  Let
$$
a_k=\inf_{x \in C^\d_n} \frac {\ln |Df^k(x)|} {k}.
$$

Analogously to
Lemma \ref {very good partial time}, we define $m_k$ as the
number of full returns
before $k$, so that $m_k$ is the biggest integer such that
$$
\sum_{i=1}^{m_k} r_n(j_i) \leq k.
$$
We define
$$
\beta_k=\sum_{\ntop {1 \leq i \leq m_k,}
{j_i \in VG(n_0,n)}} r_n(j_i),
$$
(counting the time up to $k$ spent in complete very good returns)
and
$$
i_k=k-\sum_{i=1}^{m_k} r_n(j_i).
$$
(counting the time in the incomplete return at $k$). 

Let us now consider two cases: either all iterates are part of very good
returns (that is, all $j_i$, $1
\leq i \leq m_k$ are very good and if $i_k>0$ then $j_{m_k+1}$ is also very
good), or some iterates are not part of very good returns.

{\it Case 1 (All iterates are part of very good returns).}
Since full good returns are very hyperbolic by G1 and
very good returns are good, we
just have to worry with possibly losing hyperbolicity in the incomplete
time.  To control this, we introduce the queue
$$
q_k=\inf_{x \in C^\d_n} \ln |Df^{i_k} \circ f^{k-i_k}(x)|.
$$
We have $-q_k \leq -\ln(c_{n-1}^{1/3} c_{n-1}^5)$
by Lemma \ref {lower bound} and VG,
using that the incomplete time is in the middle of a very good branch.
Let us split again in two cases: $i_k$ big or otherwise.

{\it Subcase 1a ($i_k \geq c_{n-1}^{-4/(n-1)}$).}
If the incomplete time is big, we can use
G2 to estimate the hyperbolicity of the incomplete time (which is part of a
very good return): $q_k/i_k>\lambda_{n_0}/2$.  We have
$$
a_k>\lambda_{n_0} \frac {(1+2^{n_0-n})} {2} \cdot \frac {k-i_k} {k}+
\frac {q_k} {i_k} \cdot \frac {i_k} {k}>\frac {\lambda_{n_0}} {2}.
$$

{\it Subcase 1b ($i_k<c_{n-1}^{-4/(n-1)}$).}
If the incomplete time is not big, we can not use G2 to estimate $q_k$, but
in this case $i_k$ is much less than $k$: since $k>c_{n-1}^{-4/(n-1)}$, at
least one return was completed ($m_k \geq 1$), and since it must be very
good we conclude that $k>c_{n-1}^{-1/2}$ by Remark \ref {very good time
again}, so for $n$ big
$$
a_k>\lambda_{n_0} \frac {(1+2^{n_0-n})} {2} \cdot \frac {k-i_k} {k}-
\frac {-q_k} {k}>\frac {\lambda_{n_0}} {2}.
$$

{\it Case 2 (Some iterates are not part of a very good return).}
By LC1, $m_k>c_{n-1}^{-1/30}$.
Notice that by LC2, if $m_k>c_{n-1}^{-\epsilon/5}$ then
$$
k-i_k>c_{n-1}^{-1+2\epsilon} m_k/2.
$$
So it follows that $m_k>c_{n-1}^{-1/30}$ implies that
$k>c_{n-1}^{-35/34}$ (using small $\epsilon$). 

For the incomplete time we have $-q_k \leq -\ln(c_n c_{n-1}^5) <
c_{n-1}^{-1-\epsilon}$, so $-q_k/k<c_{n-1}^{1/100}$.

Arguing as in Lemma \ref {very good partial time}, we split $k-\beta_k-i_k$
(time of full returns which are not very good) in a part relative to
returns with high time (more than $c_{n-1}^{-1-2 \epsilon}$) which we denote
$h_k$ and in a part relative to returns with low time (less than
$c_{n-1}^{-1-2\epsilon}$) which we denote $l_k$.  Using LC4 and LC5
to bound the number of returns with high time, and using LS2 to bound their
time, we get
$$
h_k<c_{n-1}^{-14} (6 \cdot 2^n) c_{n-1}^{100} m_k,
$$
and using LC1 and LC3 we have
$$
l_k<c_{n-1}^{-1-2\epsilon} (6 \cdot 2^n)
c_{n-1}^{1/60} m_k<c_{n-1}^{-79/80} m_k,
$$
provided $\epsilon$ is small enough.

Since $k>c_{n-1}^{-1+2\epsilon} m_k/2$ we have
$$
\frac {h_k+l_k} {k}<4 c_{n-1}^{1/85},
$$
provided $\epsilon$ is small enough.

Now if $i_k<c_{n-1}^{-1-2\epsilon}$ then $i_k/k<c_{n-1}^{1/80}$ (using
$\epsilon$ small), and if $i_k>c_{n-1}^{-1-2\epsilon}$ then by LC5,
$m_k \geq e^{c_{n-1}^{-\epsilon/5}}>c_{n-1}^{-n}$, so by LS2,
$i_k/k<i_k/m_k<c_{n-1}^{-14}/c_{n-1}^{-n}$.  So in both
cases $i_k/k<c_{n-1}^{1/80}$.

From our estimates on $i_k$ and on $h_k$ and $l_k$ we have
$1-(\beta_k/k)<c_{n-1}^{1/90}$.
Now very good returns are very hyperbolic, and full returns
(even not very good ones) always give derivative at least
$1$ from Lemma \ref {distortion}, so we have the estimate
$$
a_k>\lambda_{n_0} \frac {(1+2^{n-n_0})} {2} \cdot \frac {\beta_k} {k}-
\frac {-q_k} {k}>\frac {\lambda_{n_0}} {2}
$$
for $n$ big.
\end{pf}

\section{Main Theorems}

\subsection{Proof of Theorem A}


We must show that with total probability, $f$ is Collet-Eckmann.  We will
use the estimates on hyperbolicity of cool landings to show that if the
critical point always falls in a cool landing then there is uniform control
of the hyperbolicity along the critical orbit.

Let
$$
a_k=\frac {\ln |Df^k(f(0))|} {k}
$$
and $e_n=a_{v_n-1}$.

It is easy to see that if $n_0$ is big enough such that the conclusions of
both Lemmas \ref {lands cool} and \ref {cool hyperbolicity} are valid,
we obtain for $n$ large enough that
$$
e_{n+1} \geq e_n \frac {v_n-1} {v_{n+1}-1}+\frac {\lambda_{n_0}} {2} \cdot
\frac {v_{n+1}-v_n} {v_{n+1}-1}
$$
and so
\begin{equation} \label {g3}
\liminf_{n \to \infty} e_n \geq \frac {\lambda_{n_0}} {2}.
\end{equation}

Let now $v_n-1<k<v_{n+1}-1$.  Define
$q_k=\ln |Df^{k-v_n}(f^{v_n}(0))|$.

Assume first that $k \leq
v_n+c_{n-1}^{-4/(n-1)}$.  From LC1 we know that $\tau_n$ is very good, so by
LS1 we have $r_n(\tau_n)>c_{n-1}^{-1/2}$,
so $k$ is in the middle of this branch (that is, $v_n \leq k \leq
v_n+r_n(\tau_n)-1$).
Using that $|R_n(0)|>|I_n|/2^n$ (see Lemma \ref {away from the boundary}),
we get by Lemma \ref {lower bound} that
$$
-q_k<-\ln(2^{-n} c_{n-1} c_{n-1}^5)<c_{n-2}^{-1-\epsilon}.
$$
Since $v_n>c_{n-1}^{-1+\epsilon}$ (by Lemma \ref {growth of v_n}) we have
\begin{equation} \label {g1}
a_k \geq e_n \frac {v_n-1} {k}-\frac {-q_k} {k}>\left (1-\frac {1} {2^n}
\right )e_n-\frac {1} {2^n}.
\end{equation}

If $k>v_n+c_{n-1}^{-4/(n-1)}$, using Lemma \ref {cool hyperbolicity}
we get
\begin{equation} \label {g2}
a_k \geq e_n \frac {v_n-1} {k}+\frac {\lambda_{n_0}} {2} \cdot
\frac {k-v_n+1} {k}.
\end{equation}

It is clear that estimates (\ref {g3}), (\ref {g1}) and (\ref {g2})
imply that $\liminf_{k \to \infty} a_k \geq \lambda_{n_0}/2$ and so
$f$ is Collet-Eckmann.

\subsection{Proof of Theorem B} \label {proof b}

We must obtain, with total probability, upper and lower (polynomial)
bounds for the recurrence of the critical orbit.  It will be easier to
first study the recurrence with respect to iterates of return branches,
and then estimate the total time of those iterates.

\subsubsection{Recurrence in terms of return branches}

The principle of the phase analysis
is very simple: for the essentially Markov process
generated by iteration of the non-central branches of $R_n$,
most orbits (in the qs sense)
approach $0$ at a polynomial rate before falling in $I_{n+1}$.
From this we conclude,
using the phase-parameter relation, that with total probability the
same holds for the critical orbit.

\begin{lemma}

With total probability, for $n$ big enough and for
$1 \leq i \leq c^{-2}_{n-1}$,
$$
\frac {\ln |R_n^i(0)|} {\ln(c_{n-1})} <
(1+4\epsilon) \left( 1+\frac {\ln(i)} {\ln(c_{n-1}^{-1})} \right).
$$

\end{lemma}

\begin{pf}

Notice that due to torrential (and monotonic) decay of $c_n$, we can
estimate $|I_n|=c_{n-1}^{1+\delta_n}$, with $\delta_n$ decaying torrentially
fast.

From Lemma \ref {away from the boundary}, we have
$$
\frac {\ln |R_n(0)|} {\ln c_{n-1}}<\frac {\ln (2^{-n}|I_n|)} {\ln
c_{n-1}}<1+4\epsilon
$$
and the result follows for $i=1$.

For $1 \leq j \leq 2 \epsilon^{-1}$, let $X_j \subset I_n$
be a $c_{n-1}^{(1+2\epsilon)(1+j \epsilon)}$ neighborhood of $0$.
For $n$ big, we can estimate (due to the relation between $|I_n|$ and
$c_{n-1}$)
$$
\frac {|X_j|} {|I_n|}<\frac {c_{n-1}^{(1+2\epsilon)(1+j\epsilon)}}
{c_{n-1}^{1+2\epsilon}}=c_{n-1}^{j\epsilon(1+2 \epsilon)}
$$
(we of course consider $X_j$ as a union of $C^\d_n$, so that its
size is near the required size, the precision is high enough for our
purposes due to Remark \ref {approx}).

We have to make sure that the critical point does not land in some $X_j$ for
$c_{n-1}^{(1-j)\epsilon}<i \leq c_{n-1}^{-j\epsilon}$.  This requirement can
be translated on $R_n(0)$ not belonging to a certain set
$Y_j \subset I_n$ such that
$$
Y_j=\bigcup_{c_{n-1}^{(1-j)\epsilon} \leq |\d|<c_{n-1}^{-j\epsilon}}
(R_n^\d)^{-1}(X_j).
$$

By Lemma \ref {away from the boundary},
it is clear that no $X_j$ intersects $I^{\tau_n}_n$, so we easily get
$$
p_\g(Y_j|I^{\tau_n}_n) \leq
c_{n-1}^{-j\epsilon} c_{n-1}^{(1+\epsilon) j\epsilon} \leq
c_{n-1}^{\epsilon^2}
$$
and
$$
p_\g(\bigcup_{j=1}^{2 \epsilon^{-1}} Y_j|I^{\tau_n}_n) <
2\epsilon^{-1}c_{n-1}^{\epsilon^2}.
$$

Applying PhPa1, the probability that for some
$1 \leq j \leq 2\epsilon^{-1}$ and $c^{(1-j)\epsilon}_{n-1} < i \leq
c^{-j\epsilon}_{n-1}$ we have $|R^i_n(0)|<c^{(1+2
\epsilon)(1+j\epsilon)}_{n-1}$ is at most
$2\epsilon^{-1}c_{n-1}^{\epsilon^2}$, which is summable.
In particular, with total probability, for $j$ and $i$ as above, we have
for $n$ big enough
\begin{align*}
\frac {\ln |R_n^i(0)|} {\ln(c_{n-1})}
&\leq (1+2\epsilon)(1+j\epsilon)\\
&<(1+4\epsilon)(1+(j-1)\epsilon) <
(1+4\epsilon) \left( 1+\frac {\ln(i)} {\ln(c_{n-1}^{-1})} \right).
\end{align*}
\end{pf}

\begin{lemma}

With total probability, for $n$ big enough and for
$c_{n-1}^{-1-\epsilon}<i \leq s_n$,
$$
\frac {\ln |R_n^i(0)|} {\ln(c_{n-1})} <
(1+4\epsilon) \left( 1+\frac {\ln(i)} {\ln(c_{n-1}^{-1})} \right).
$$

\end{lemma}

\begin{pf}

The argument is the same as for the previous lemma, but the decomposition
has a slightly different geometry.  Let
$$
x_j=c_{n-1}^{(1+2\epsilon)(1+(1+\epsilon)^{j+1})},
$$
so that
$$
\frac {x_j} {|I_n|}<\frac {c_{n-1}^{(1+2\epsilon)(1+(1+\epsilon)^{j+1})}}
{c_{n-1}^{1+\epsilon}}<c_{n-1}^{(1+2\epsilon)(1+\epsilon)^{j+1}}.
$$
Let $K$ be biggest with $x_K>c_n^{1-\epsilon}$.
For $0 \leq j \leq K$, let $X_j \subset I_n$
be a $x_j$ neighborhood of $0$
(approximated as union of $C^\d_n$, notice that $x_j>c_n^{1-\epsilon} \gg
|I_{n+1}|$ for $0 \leq j \leq K$, so the approximation
is good enough for our purposes due to Remark \ref {approx}).
Let $Y_j \subset I_n$ be such that
$$
Y_j=\bigcup_{c_{n-1}^{-(1+\epsilon)^j} \leq |\d|<
c_{n-1}^{-(1+\epsilon)^{j+1}}} (R_n^\d)^{-1}(X_j).
$$

By Lemma \ref {away from the boundary}, it is clear that no
$X_j$ intersects $I^{\tau_n}_n$, so we easily get
$$
p_\g(Y_j|I^{\tau_n}_n) \leq
c_{n-1}^{-(1+\epsilon)^{j+1}}c_{n-1}^{(1+\epsilon)^{j+2}} \leq
c_{n-1}^{\epsilon(1+j\epsilon)}
$$
and
$$
p_\g(\bigcup_{j=0}^K Y_j|I^{\tau_n}_n)<\sum_{j=0}^\infty
c_{n-1}^{\epsilon(1+j\epsilon)}=\frac {c_{n-1}^\epsilon}
{1-c_{n-1}^{\epsilon^2}}<c_{n-1}^{\epsilon/2}.
$$

Applying PhPa1, the probability that for some
$0 \leq j \leq K$ and
$$
c^{-(1+\epsilon)^j}_{n-1} < i \leq
c^{-(1+\epsilon)^{j+1}}
$$
we have
$$
|R^i_n(0)|<c^{(1+2 \epsilon)(1+(1+\epsilon)^{j+1})}_{n-1}
$$
is at most
$c_{n-1}^{\epsilon/2}$, which is summable.
In particular, with total probability, for $j$ and $i$ as above, we have
\begin{align*}
\frac {\ln |R_n^i(0)|} {\ln(c_{n-1})}
&<(1+2\epsilon)(1+(1+\epsilon)^{j+1})\\
&<(1+4\epsilon)(1+(1+\epsilon)^j)<
(1+4\epsilon) \left( 1+\frac {\ln(i)} {\ln(c_{n-1}^{-1})} \right).
\end{align*}

This covers the range $c_{n-1}^{-1} < i \leq c_{n-1}^{-(1+\epsilon)^{K+1}}$. 
For $c_{n-1}^{-(1+\epsilon)^{K+1}} < i \leq s_n$, notice that
$R_n^i(0) \notin I_{n+1}$, so
\begin{align*}
\frac {\ln |R_n^i(0)|} {\ln c_{n-1}}
&<\frac {\ln (|I_{n+1}|/2)} {\ln c_{n-1}}\\
&<\frac {1+4\epsilon} {1+2\epsilon} \cdot \frac {\ln c_n^{1-\epsilon}} {\ln
c_{n-1}}\\ 
&\leq \frac {1+4\epsilon} {1+2\epsilon} \cdot \frac {\ln x_{K+1}} {\ln
c_{n-1}} &&\text {(by definition of $K$)}\\
&\leq (1+4\epsilon)(1+(1+\epsilon)^{K+1})\\
&\leq (1+4\epsilon) \left( 1+\frac {\ln(i)} {\ln(c_{n-1}^{-1})} \right).
\end{align*}
\end{pf}

Both cases are summarized below:

\begin{cor} \label {distance}

With total probability, for $n$ big enough and for
$1 \leq i \leq s_n$,
$$
\frac {\ln |R_n^i(0)|} {\ln(c_{n-1})} <
(1+4\epsilon) \left( 1+\frac {\ln(i)} {\ln(c_{n-1}^{-1})} \right).
$$

\end{cor}

\subsubsection{Total time of full returns}

We must now relate the return times in terms of $R_n$ to the return
times in terms of $f$.

For $1 \leq i \leq s_n$, let $k_i$ be
such that $R_n^i(0)=f^{k_i}(0)$.

\begin{lemma} \label {rec1}

With total probability, for $n$ big enough and for
$c_{n-1}^{-\epsilon}<i \leq s_n$,
$$
\frac {\ln(k_i)} {\ln(c_{n-1}^{-1})}>(1-3\epsilon) \left (1+\frac {\ln(i)}
{\ln(c_{n-1}^{-1})} \right ).
$$
\end{lemma}

\begin{pf}

By Lemma \ref {lands cool}, $R_n(0)$ belongs to a cool landing, so
using LC2 (which allows to estimate the average of return times over a large
initial segment of cool landings) we get
$$
\frac {k_i} {i-1}>c_{n-1}^{-1+3\epsilon}.
$$
This immediately gives
$$
\frac {\ln(k_i)} {\ln(c_{n-1}^{-1})}>(1-3\epsilon) + \frac {\ln(i-1)}
{\ln(c_{n-1}^{-1})}>(1-3\epsilon) \left (1+\frac {\ln(i)}
{\ln(c_{n-1}^{-1})} \right )
$$
\end{pf}

Using that $v_n>c_{n-1}^{-1+\epsilon}$
(from Corollary \ref {growth of v_n})
and that $k_i \geq v_n$ we get for $1 \leq i \leq c_{n-1}^{-\epsilon}$
$$
\frac {\ln(k_i)} {\ln(c_{n-1}^{-1})} \geq
\frac {\ln(v_n)} {\ln(c_{n-1}^{-1})}>
\frac {\ln(c_{n-1}^{-1+\epsilon})} {\ln(c_{n-1}^{-1})}>
(1-3\epsilon) (1+\epsilon) \geq
(1-3\epsilon) \left( 1+\frac {\ln(i)} {\ln(c_{n-1}^{-1})} \right).
$$
Together with Lemma \ref {rec1}, this gives

\begin{cor} \label {high time}

With total probability, for $n$ big enough and for
$1 \leq i \leq s_n$,
$$
\frac {\ln(k_i)} {\ln(c_{n-1}^{-1})} >
(1-3\epsilon) \left( 1+\frac {\ln(i)} {\ln(c_{n-1}^{-1})} \right).
$$

\end{cor}

\subsubsection{Upper and lower bounds}

Notice that
$|R_n(0)|=|f^{v_n}(0)| \leq c_{n-1}$, so using Lemma \ref {growth of
v_n} we get
$$
\limsup_{n \to \infty} \frac {-\ln |f^n(0)|}{\ln (n)} \geq
\limsup_{n \to \infty} \frac {-\ln |f^{v_n}(0)|}{\ln (v_n)} \geq
\limsup_{n \to \infty} \frac {-\ln (c_{n-1})}{\ln (v_n)} \geq 1.
$$

Let now $v_n \leq k<v_{n+1}$.  If $|f^k(0)|<k^{-1-10\epsilon}$ then Lemma
\ref {growth of v_n} implies that
$f^k(0) \in I_n$ and so $k=k_i$ for some $i$.
It follows from Corollaries \ref {distance} and \ref {high time} that
$$
|f^{k_i}(0)|>k_i^{-1-10 \epsilon}.
$$

Varying $\epsilon$ we get
$$
\limsup_{n \to \infty} \frac {-\ln |f^n(0)|}{\ln (n)} \leq 1.
$$

\appendix

\section{Sketch of the proof of the Phase-Parameter relation}

The proof of the Phase-Parameter relation uses ideas from complex analysis. 
We will provide a sketch of the proof assuming familiarity with the
work \cite {parapuzzle}.  For a more general result (with all details fully
worked out), see \cite {AM4}.

Given a simple map $f$, one can define (as in \S 3 of \cite {parapuzzle})
a sequence of {\it holomorphic
families of generalized quadratic-like maps} $R_i$, $i \geq 1$,
related by {\it generalized renormalization}.  To fix notation, the
parameter space of those families will be denoted $\Lambda_i[f]$, so for
each $g \in \Lambda_i[f]$ the family defines
a {\it generalized quadratic-like
map} $R_i[g]:U^j_i[g] \to U_i[g]$.  Moreover, the family $R_i$ is {\it
equipped} (with a holomorphic motion $h_i$ of the $U^j_i$ and $U_i$) and
{\it proper}.  The following properties of this sequence of families will be
important for us:

(1)\, $\Lambda_i[f] \cap \R=J_i[f]$, and
for $g \in J_i[f]$, $U_i[g] \cap \R=I_i[g]$ and $U^j_i[g] \cap
\R=I^j_i[g]$.

(2)\, For $g \in J_i[f]$, the map $R_i[g]:\cup U^j_i[g] \to U_i[g]$
is an extension of the {\it real first return map}
$R_i[g]:I^j_i[g] \to I_i[g]$ defined before.
For $j \neq 0$ (respectively for $j=0$),
$R_i[g]:U^j_i[g] \to U_i[g]$ is a homeomorphism (respectively a double
covering).

(3)\, The modulus of $U_i[f] \setminus \overline {U_{i+1}[f]}$ grows at
least linearly in $i$ (\cite {GS2}, \cite {puzzle}).

(4)\, The modulus of $\Lambda_i[f] \setminus \overline {\Lambda_{i+1}[f]}$
grows at least linearly in $i$ (\cite {parapuzzle}).

Define $\tau_i$ as before by $R_i[f](0) \in I^{\tau_i}_i[f]$.
Let $\Lambda^j_i[f]$ denote the set of
$g \in \Lambda_i[f]$ such that $R_i[g](0) \in
U^j_i[g]$.  We have:

(5)\, $\Lambda^j_i[f] \cap \R=J^j_i[f]$ and in particular
$\Lambda^{\tau_i}_i[f] \cap \R=J^{\tau_i}_i[f]$.

By Lemma 4.8 of \cite {parapuzzle}, the holomorphic motion $h_i$ of
$U_i$, $U^j_i$ (corresponding to $R_i$) has uniformly bounded dilatation
(independently of $i$) {\it when restricted to
$U_i \setminus \overline {U^0_i}$}.
Item (3) above implies that for $i$ big, there is an annulus
of big modulus (linear growth in $i$)
contained in $U_i[f] \setminus \overline {(U^0_i[f] \cup
U^{\tau_i}_i[f])}$ and going around $U^{\tau_i}_i[f]$.  By transverse
quasiconformality of holomorphic motions
(Corollary 2.1 of \cite {parapuzzle}), and the
$\lambda$-Lemma of \cite {MSS},
this estimate can be transferred to the parameter space, so we get:

(6)\, The modulus of $\Lambda_i[f] \setminus \overline
{\Lambda^{\tau_i}_i[f]}$ grows at least linearly in $i$.

For each $g \in \Lambda_i[f]$, denote by $L_i[g]$ the first landing map to
$U^0_i[g]$ obtained by iteration of non-central branches of $R_i[g]$.  By
item (2) above, we have:



(7)\, For $g \in J_i[f]$, the domain of $L_i[g]$ is a union $\cup_{\d \in
\Omega} W^\d_i[g]$ of disks such that $W^\d_i[g] \cap \R=C^\d_i[g]$ and
$L_i[g]$ extends the {\it real first landing map}
$L_i[g]:\cup_{\d \in \Omega} C^\d_i[g] \to I^0_i[g]$ defined before.

The family $L_i$ is also equipped with a holomorphic motion
$\hat h_i$ of $U_i$
and the $W^\d_i$ (see \S 3.5 of \cite {parapuzzle}).
Define $\Gamma^\d_i[f]$ as the set of $g
\in \Lambda_i[f]$ such that $R_i[g](0) \in W^\d_i[g]$.
The $\lambda$-Lemma and (6) imply:

(8)\, For $g \in J^{\tau_i}_i[f]$, there exists a
real-symmetric qc map of $\C$, whose dilation goes to $1$ as $i$
grows, taking $U_i[f]$ to
$U_i[g]$, and taking any $W^\d_i[f]$ to $W^\d_i[g]$.

Items (7) and (8) prove PhPh1 in the Phase-Parameter relation.

Item (6) and transverse quasiconformality of holomorphic motions imply:

(9)\, There is a real-symmetric qc map of $\C$, whose
dilatation goes to $1$ as $i$ grows,
taking $U^{\tau_i}_i[f]$ to $\Lambda^{\tau_i}_i[f]$, and taking any
$W^\d_i[f]$ {\it contained in} $U^{\tau_i}_i[f]$ to $\Gamma^\d_i[f]$.

Items (5), (7) and (9) prove PhPa1 in the Phase-Parameter relation.

\bigskip

Notice that for any $g \in \Lambda_i[f]$, and for every $\d$ the map
$L_i[g]:W^\d_i[g] \to U^0_i[g]$ extends to a holomorphic diffeomorphism
$R^\d_i[g]:U^\d_i[g] \to U_i[g]$.  It is easy to see that $\hat h_i$
(as defined in \S 3.5 of \cite {parapuzzle})
is also a holomorphic motion of the $U^\d_i$.

Define $\tilde \Lambda_{i+1}[f]$ as
the set of $g$ such that $R_i[g](0) \in U^{\d_i}_i[g]$, where $\d_i$ is
chosen such that $R_i[f](0) \in C^{\d_i}_i[f]$.  It follows that
$\Lambda_{i+1}[f]=\Gamma^{\d_i}_i[f] \subset
\tilde \Lambda_{i+1}[f] \subset
\Lambda^{\tau_i}_i[f]$.  By (3), the modulus of $U^{\d_i}_i[f] \setminus
\overline {W^{\d_i}_i[f]}$ grows at least
linearly in $i$.  By (6) and transverse
quasiconformality of holomorphic motions, this implies:

(10)\, The modulus of
$\tilde \Lambda_i[f] \setminus \overline {\Lambda_i[f]}$
grows at least linearly in $i$.

For $g \in \Lambda_i[f]$, the map
$R_i[g]:U^0_i[g] \to U_i[g]$ extends to a bigger domain
$\tilde U_{i+1}[g]=(R_{i-1}[g]|_{U_i[g]})^{-1}(U^{\d_{i-1}}_{i-1}[g])$,
as a double covering map onto
$U_{i-1}[g]$ (notice that $U^0_i[g] \subset \tilde
U_{i+1}[g] \subset U_i[g]$).  It follows:

(11)\, If $g \in J_i[f]$ then $\tilde U_{i+1}[g] \cap \R=\tilde
I_{i+1}[g]$.



The holomorphic motion $\hat h_{i-1}$ (corresponding to $L_{i-1}$)
naturally lifts to a holomorphic motion $\tilde h_i$ of
$U_i$, $\tilde U_{i+1}$ and all $U^j_i$ {\it not contained
in} $\tilde U_{i+1}$, which is defined (in principle) over $\Lambda_i[f]$,
but extends to a holomorphic motion defined over $\tilde \Lambda_i[f]$.



Item (10) and yet another application of the $\lambda$-Lemma imply:

(12)\, For $g \in J_i[f]$, there exists a
real-symmetric qc map of $\C$, whose dilatation goes to $1$ as
$i$ grows, taking $U_i[f]$ to
$U_i[g]$, and taking any $U^j_i[f]$ {\it not
contained in} $\tilde U_{i+1}[f]$ to $U^j_i[g]$.

Items (2), (11) and (12) prove PhPh2 in the Phase-Parameter relation.

Item (10) and transverse quasiconformality of holomorphic motions imply:

(13)\, There is a real-symmetric qc map of $\C$, whose
dilatation goes to $1$ as $i$ grows,
taking $U_i[f]$ to $\Lambda_i[f]$, and taking any
$U^j_i[f]$ {\it not contained in} $\tilde U_{i+1}[f]$ to
$\Lambda^j_i[f]$.

Items (2), (11) and (13) prove PhPa2 in the Phase-Parameter relation.
All items of the Phase-Parameter relation are proved.

\end{document}